\documentclass{article}
\usepackage{amsthm}    
\usepackage{amsmath,amssymb,amsfonts} 
\usepackage{graphicx} 
\usepackage{enumitem} 
\usepackage{hyperref} 
\usepackage{cleveref} 
\usepackage{algorithm}    
\usepackage{algpseudocode} 
\usepackage{float} 
\usepackage{graphicx}
\usepackage{subcaption}
\usepackage{cleveref} 
\usepackage{booktabs,multirow}

\hypersetup{colorlinks=true, linkcolor=blue, citecolor=blue}
\crefname{lemma}{Lemma}{Lemmas}
\crefname{theorem}{Theorem}{Theorems}
\crefname{definition}{Definition}{Definitions}
\crefname{equation}{equation}{equations}
\newcommand{\norm}[1]{\|#1\|} 
\newcommand{\supp}{\operatorname{supp}} 
\newcommand{\argmin}{\operatorname{argmin}} 
\newcommand{\dom}{\operatorname{dom}} 
\newcommand{\intr}{\operatorname{int}} 

\newtheorem{definition}{Definition}[section] 
\newtheorem{lemma}{Lemma}[section] 
\newtheorem{theorem}{Theorem}[section] 
\newtheorem{example}{Example}[section]
\newtheorem{remark}{Remark}[section]
\newtheorem{assumption}{Assumption}[section]

\title{BGPNP: A Non-Euclidean Gradient Projection Newton Pursuit Method for Sparsity-Constrained  Optimization}

\author{Jiawei Chen\footnote{ School of Mathematics and Statistics, Southwest
University, Chongqing, 400715, China, E-mail: j.w.chen713@163.com.},\quad Huilin Han\footnote{ School of Mathematics and Statistics, Southwest
University, Chongqing, 400715, China,  E-mail: hhllhh113@163.com.},\quad Yibing Lv\footnote{School of Information and Mathematics, Yangtze University, Jingzhou 434023, China,  E-mail: Yibinglv@yangtzeu.edu.cn.},
\quad
Jen-Chih Yao\footnote{Research Center for Interneural Computing, China Medical University, Taichung 40447, Taiwan, E-mail: yaojc@mail.cmu.edu.tw.}}

\date{}

\begin{document}
		\maketitle
\begin{abstract}
 A non-Euclidean gradient projection Newton pursuit method is proposed for solving sparsity constrained optimization problem via the Bregman distance. The method is also denoted as  BGPNP. The method firstly constructs a gradient projection step based on the Bregman distance to ensure that the iterates satisfy the sparsity constraint. Then, when the support set becomes stable or the gradient norm is sufficiently small, a Newton pursuit step is introduced to improve local convergence efficiency.
 The BGPNP is proved to possess a sufficient descent property and the non-increasing monotonicity of the objective function values.
  The fully convergence  and linear convergence rate  of the BGPNP  are  established under appropriate conditions.  Moreover,
  the locally
quadratical convergence rate of the BGPNP is also obtained.
  This method provides an effective algorithmic framework for sparse optimization in non-Euclidean settings. Finally, the
  numerical experiments on sparsity-constrained compressed sensing problems validate the competitive performance of BGPNP compared with existing leading solvers.
  \end{abstract}
		

\section{Introduction}
	\label{sec:introduction}
In many real-world problems, the available data or decision variables are often high-dimensional, while only a small number of variables contribute significantly to the underlying structure or decision-making process. This characterization, known as sparsity, reflects the fact that meaningful information is frequently contained in a small subset of variables. Exploiting such sparse structures can reduce computational complexity, improve model interpretability, and enhance the efficiency of practical applications. A variety of applications can be naturally formulated as sparse optimization problems. Just because of this,  sparse optimization has attracted considerable attention in optimization theory,
statistics, economics, communication, AI and data science, such  as compressed sensing  \cite{eldar2012compressed,foucart2013invitation}, machine learning \cite{guyon2003introduction,tibshirani1996regression}, signal processing \cite{GarmroudiEH}, phase retrieval \cite{DaiLY2026,FanYXL},  high-dimensional statistical inference \cite{KohKimBoyd,MingYX,Moosaei} and financial portfolio optimization \cite{bertsimas2009algorithm,chang2000heuristics,markowitz1952portfolio}.
For instance, in compressed sensing, sparse recovery aims to reconstruct signals from a limited number of measurements by exploiting their intrinsic sparse representations  \cite{bruckstein2009sparse,candes2006stable,donoho2006compressed}. In machine learning and statistical modeling, sparse optimization provides an effective approach for feature selection by identifying the most relevant variables from large-scale datasets \cite{guyon2003introduction,KohKimBoyd,MingYX,tibshirani1996regression}. In financial portfolio optimization, sparsity constraints allow investors to select a limited number of assets, thereby reducing transaction costs and improving portfolio interpretability \cite{bertsimas2009algorithm,chang2000heuristics,markowitz1952portfolio}.

Sparse optimization is to find a sparse solution of  optimization problems with sparse constraints.
A basic model for sparse optimization is the $\ell_0$-norm constrained optimization problem:
\begin{equation}\label{qu:1}
	\begin{split}
		&\underset{x}{\operatorname{min}} \, f(x)\\
		&\text{s.t.}\, \|x\|_{0} \leq s,
	\end{split}
\end{equation}
where $\|x\|_{0}$ denotes the number of nonzero entries of $x$, and $s$ represents the prescribed sparsity level. The $\ell_0$-norm constraint directly characterizes the sparsity of solutions and provides an intuitive formulation for variable selection and sparse recovery problems.
However, sparse optimization problem is NP-hard even for linear problems due
to the discontinuity and nonconvexity of the sparse constraint and the combinatorial properties of sparse feasible region.
Recent years, the theory, such as optimality conditions, duality, saddle point theorem and constraint qualifications, of sparsity-constrained optimization (SCO) have been widely investigated by leveraging the variational properties of the sparse constraints; see e.g., \cite{BeckE2013,ChenSOL24,ChenWLY2026,KanSong2024,PanLuoXiu2017,panXiuF,ZhaoLuoLQX} and the references therein.
It laid the theoretical foundation for development of numerical algorithms of sparse optimization problems.

%

Due to the nonconvex and discontinuous nature of the $\ell_0$-norm and the combinatorial properties of sparse feasible region, it is generally challenging to solve sparsity constrained optimization problems. A large class of efficient numerical methods is based on the hard-thresholding principle, which enforces sparsity by retaining only the largest components of the iterates. The representative algorithms include gradient support pursuit (GraSP) \cite{bahmani2013greedy}, iterative hard thresholding (IHT)\cite{blumensath2009iht}, normalized iterative hard thresholding (NIHT)\cite{blumensath2010normalized}, subspace pursuit (SP)\cite{dai2009subspace}, hard thresholding pursuit (HTP)\cite{foucart2011htp},  compressive sampling matching pursuit (CoSaMP)\cite{needell2009cosamp}, and gradient hard thresholding pursuit (GraHTP)\cite{yuan2018gthp}. These methods are computationally efficient and easy to implement, but their local convergence speed is often limited.
{\it An interesting issue is how to accelerate the local convergence behavior of these algorithms}.

A natural way is to make full use of the properties, such as second-order information, of the underlying functions and the structure of sparse optimization.
 In particular, Newton-type hard-thresholding methods exploit curvature information restricted to the active support set, leading to significantly faster local convergence. Wu, Pan and Yang \cite{WuPanYang2023} proposed a hybrid of the proximal gradient method and the subspace regularized Newton method for $\ell_{q}$-norm composite optimization problem, and obtained the superlinear convergence rate for the proposed method under an additional local error bound condition
and Kurdyka-\L{}ojasiewicz property of the cost function with exponent $1/2$.  Zhou, Xiu and Qi \cite{zhou2021global} proposed a Newton hard-thresholding pursuit method with global and quadratic convergence guarantees for SCO. Zhou et al \cite{zhou2022gradient} developed a gradient projection Newton pursuit method (GPNP) for SCO, which combines first-order sparse projection steps with selectively activated Newton steps, achieving a balance between computational efficiency and local acceleration. These works show that the convergence behavior of first-order methods can be improved by using the second-order information. However, the full convergence and locally quadratical convergence of the proposed methods \cite{zhou2022gradient,zhou2021global} are established under the restricted strong smoothness of the objective function (i.e., the gradient Lipschitz continuity of objective function) in the sense of Euclidean norm. This condition may confine the applicability of the proposed methods \cite{zhou2022gradient,zhou2021global}.
In particular, for objective functions whose gradients do not satisfy conventional Lipschitz continuity assumptions but exhibit favorable smoothness properties relative to an appropriate reference function, the standard quadratic upper bound may not provide a suitable description of their behavior. Therefore, it is  worthwhile to extend the frameworks of \cite{WuPanYang2023,zhou2022gradient,zhou2021global} beyond Euclidean distance and develop a unified approach to sparse projection and Newton pursuit under more general distance structures.

Bregman distances, introduced by Bregman \cite{bregman1967relaxation}, provide a flexible non-Euclidean framework induced by a convex kernel function.
In recent years, Bregman distances have been widely used in mirror descent, proximal algorithms, and non-Euclidean gradient methods \cite{bauschke2019linear,bauschke2017descent,beck2003mirror}.
Bauschke, Bolte, and Teboulle \cite{bauschke2017descent} established a descent lemma beyond conventional Lipschitz gradient continuity by introducing relative smoothness, thereby providing a fundamental theoretical basis for Bregman proximal gradient methods.
In particular, Bauschke et al. \cite{bauschke2019linear} established linear convergence of Bregman gradient methods under the gradient dominance conditions, without requiring classical assumptions such as strong convexity or Lipschitz gradient continuity.
Subsequently, Dragomir et al. \cite{Dragomir2022} investigated the optimal complexity and certification of Bregman first-order methods, further highlighting the role of non-Euclidean geometry in first-order optimization theory. For nonconvex optimization, Ahookhosh, Themelis, and Patrinos \cite{Ahookhosh2021} proposed a Bregman forward-backward linesearch algorithm and established superlinear convergence under suitable assumptions. Latafat et al. \cite{Latafat2022} developed the Bregman Finito/MISO method for nonconvex regularized finite-sum minimization without requiring Lipschitz continuity of the gradient.
Gao et al. \cite{Gao2023} proposed an alternating structure-adapted Bregman proximal gradient method and its inertial variant for constrained nonconvex nonsmooth optimization. Zhang et al. \cite{Zhang2021} established linear convergence of a proximal-like incremental aggregated gradient method under Bregman distance growth conditions, demonstrating the potential of non-Euclidean distance structures in convergence rate analysis. More recently, Yang and Toh \cite{Yang2026} investigated inexact Bregman proximal gradient methods and their inertial variants, while Elshiaty and Petra \cite{Elshiaty2026} developed a multilevel Bregman proximal gradient descent framework.
These studies demonstrate that Bregman distances not only extend the applicability of conventional gradient-based methods but also provide useful tools for designing optimization algorithms with more general framework.

For sparsity-constrained optimization problems, the sparsity-constrained feasible region is nonconvex, so the corresponding Bregman projection does not necessarily inherit the standard properties of projections onto convex sets and may be set-valued.
 Moreover, since Bregman distances are generally asymmetric and fail to satisfy the triangle inequality, the sufficient descent and convergence results established for the GPNP method cannot be directly extended to the non-Euclidean setting. Therefore, it remains an open and challenging issue to effectively integrate Bregman distances into the Newton pursuit framework for sparse optimization.


Motivated by the above works, this paper investigates sparsity constrained optimization and proposes a non-Euclidean Gradient Projection Newton Pursuit algorithm in the framework of Bregman distances, called BGPNP. The proposed method firstly constructs a gradient projection step based on the Bregman distance to ensure that the iterates satisfy the sparsity constraint. Then, when the support set becomes stable or the gradient norm is sufficiently small, a Newton pursuit step is introduced to improve local convergence efficiency.
Compared with classical hard-thresholding methods, the proposed algorithm incorporates non-Euclidean metric induced by Bregman distances. In contrast to the existing Bregman first-order  methods, it further integrates Newton-type acceleration to enhance local convergence behavior. In particular, when the kernel function is chosen as the quadratic function $\frac{1}{2}\|\cdot\|^{2}$, the proposed framework recovers the GPNP method as a special case.
 The theoretical convergence results, such as the sufficient descent property, the non-increasing monotonicity of the objective function values, the fully convergence  and linear convergence, of BGPNP are established under suitable conditions.
 The locally
quadratical convergence rate of the BGPNP is also derived under the Legendre strong convexity, locally restricted Hessian Lipschitz continuity  and relative smoothness conditions.  The results presented in this paper are original even if the BGPNP reduces to GPNP.
Numerical experiments on sparsity-constrained compressed sensing problems also verify the competitive performance of BGPNP compared with existing leading solvers.




The main contributions of this paper are summarized as follows:
\begin{enumerate}
	\item[\rm(i)] A Bregman gradient projection Newton pursuit (BGPNP) algorithm is proposed for solving sparse constrained optimization problem \eqref{qu:1}. The BGPNP extends the classical Euclidean GPNP algorithm \cite{zhou2022gradient} to a non-Euclidean setting by incorporating Bregman sparse projection into the GPNP framework.
	
\item[\rm(ii)]
The BGPNP is proved to satisfy sufficient descent property and converge to stationary points under relative smoothness and strong convexity assumptions. The relative smoothness is weaker than the gradient Lipschitz continuity. Furthermore, we establish the local linear convergence and local quadratic convergence of BGPNP. In particular, the local linear convergence rates of BGPNP are original even if the Legendre function $h(x)=\frac{\|x\|^2}{2}$.

\item[\rm(iii)]
 We introduce a nonquadratic Bregman kernel with an explicit inverse gradient mapping in numerical experiments. Numerical experiments demonstrate that appropriately chosen Bregman kernels can improve computational efficiency and recovery accuracy in certain high-dimensional and noisy sparse recovery problems.
\end{enumerate}

The remainder of this paper is organized as follows. Section 2 introduces the preliminaries of Bregman distance and related assumptions. Section 3 presents the proposed Bregman Gradient Projection Newton Pursuit algorithm and establishes its convergence properties. Section 4 analyzes the convergence behavior of the Newton pursuit step. Section 5 reports numerical experiments to demonstrate the effectiveness of the proposed method.

	\section{Preliminaries}
	\label{sec:preliminaries}
	
	Throughout the paper, let $\mathbb{R}^n$ be a $n$-dimensional Euclidean space with inner product $\langle \cdot , \cdot \rangle$, $f:\mathbb{R}^n\rightarrow \mathbb{R}$ be a twice continuously differentiable function and bounded from below, $h:\mathbb{R}^n\rightarrow \mathbb{R}\cup\{+\infty\} $ be a proper, lower semicontinuous and convex function. $C^2$ denotes the class of twice continuously differentiable functions. The domain of $h$ is denoted by $\dom h:= \{ x\in \mathbb{R}^n\,: \, h(x) < +\infty \}$ and its interior is denoted by $\intr \dom h$.
 For  $x\in\mathbb{R}^n$, $\operatorname{supp}(x):=\{i\, :\,x_i\neq 0, i\in\{1,2,\ldots,n\}\}$ denotes its support set. For an index set $\Gamma\subseteq\{1,2,\ldots,n\}$, $x_\Gamma$ denotes the subvector of $x$ restricted to $\Gamma$, the cardinal number of $J$ is
  denoted by $|J|$,
  $\bar{\Gamma}=\{1,2,\ldots,n\}\setminus \Gamma $   and $\nabla_\Gamma f(x)$,
	\(\nabla^2_{\Gamma\Gamma}f(x)\) denote the corresponding restricted gradient and restricted Hessian.
	Let $\|\cdot\|$ be the Euclidean norm for a vector and the spectral norm for a matrix and set
		$S:=\{x\in\mathbb{R}^n:\|x\|_0\le s\}$.

	We next recall several notions and facts from Bregman geometry.
	
	\subsection{Bregman Distance}
	\begin{definition}[Legendre function]
		A proper, lower semicontinuous, convex function \(h\) is called a Legendre function iff it is essentially smooth and essentially strictly convex on \(\operatorname{int}\operatorname{dom}h\)$\neq \emptyset$.
	\end{definition}
	
	\begin{definition}[Bregman Distance]
	\label{def:BD}
		The Bregman distance associated with a Legendre function $h$ is defined by:
	\[
	D_{h}(x, y):=h(x)-h(y)-\langle\nabla h(y), x-y\rangle, \quad \forall x \in \dom h, y \in \intr \dom h.
	\]
    \end{definition}

	A key property of Bregman distance is
	\[
	D_{h}(x, y)=0 \Leftrightarrow x=y, \quad \forall(x, y) \in \dom h \times \intr \dom h.
	\]
	
	The following properties of Bregman distance will be useful:
	\begin{align*}
		D_{h}(x, y)&=D_{h}(x, z)+D_{h}(z, y)+\langle\nabla h(z)-\nabla h(y), x-z\rangle, \\
		D_{h}(x, y)&=D_{h^{*}}(\nabla h(y), \nabla h(x)), \quad \forall x, y \in \intr \dom h, \\
		\left(\forall h_{1}, h_{2} \in \mathcal{C}^{1}\right)&\left(\forall \beta_{1}, \beta_{2} \in \mathbb{R}\right)\,\, D_{\beta_{1} h_{1}+\beta_{2} h_{2}}=\beta_{1} D_{h_{1}}+\beta_{2} D_{h_{2}}.
	\end{align*}

		For any $z\in \operatorname{int}\operatorname{dom}h$, the Bregman sparse projection onto $S$ is defined by
	\begin{align*}
	\Pi_S(z)=\arg\min_{u\in S} D_h(u,z).
	\end{align*}
	
In particular, if the kernel function $h(x)=\frac{1}{2} \|x\|^2$, this projection $	\Pi_S(z)$ reduces to the classical hard-thresholding projection.


	\subsection{Key Definitions and Lemma}
	In this subsection, the proper, lower semicontinuous and convex function $h:\mathbb{R}^n\rightarrow \mathbb{R}\cup\{+\infty\} $ is always assumed to be a Legendre function.

	\begin{definition}[LC condition]\cite{bauschke2019linear}
	The pair of functions $(f, h)$ is said to satisfy the LC condition iff there exists $L>0$ such that $Lh-f$ is convex on  $\intr \dom h$.
	\end{definition}

	\begin{lemma}[Descent lemma]{\rm \cite{bauschke2017descent}}\label{lem:1}
		For all $(x, y) \in \intr \dom h \times \intr \dom h$, the LC-condition for the pair of functions $(f,h)$ is equivalent to
		\[
		f(x) \leq f(y)+\langle\nabla f(y), x-y\rangle+L D_{h}(x, y),
		\]
		i.e.,
		   $ D_{f}(x, y) \leq L D_{h}(x, y)$.
	\end{lemma}

	\begin{definition}[Legendre strongly convex]{\rm \cite{bauschke2019linear}}
		\label{def:1}
		A function $f:\mathbb{R}^n\to \mathbb{R}$ is Legendre strongly convex with respect to a Legendre function $h$  iff there exists  $l > 0$ such that for all $x, y \in \intr \dom h$,
		\begin{equation*}\label{lsc}
			f(x) \geq f(y)+\langle\nabla f(y), x-y\rangle+l D_{h}(x, y),
		\end{equation*}			
		i.e.,
		\begin{equation*}\label{lsc1}
		D_{f}(x, y) \geq lD_{h}(x, y).
	\end{equation*}	
		 When $l=0$, $f$ is a convex function.
	\end{definition}
	
	\begin{definition}
		A function $f:\mathbb{R}^n\to \mathbb{R}$ is said to be coercive iff $f(x)\to +\infty$ whenever $\|x\|\to\infty$.
		\end{definition}

%
%
	
	\begin{definition}[Restricted Hessian Lipschitz continuous]{\rm \cite{zhou2022gradient}}
		\label{def:linear_reg}
		A function $f$ is said to be $s$-Restricted Hessian Lipschitz continuous ($s$-RHLC) iff there exists a constant $M_{r} > 0$ such that for any index set $\Gamma$ with $|\Gamma| \leq s$,
		\begin{align*}
		\|\nabla^2_{\Gamma\Gamma}f(x) - \nabla^2_{\Gamma\Gamma}f(y)\| \leq M_{r} \|x-y\|,\,\, \forall\, x, y \in \mathbb{R}^n.
		\end{align*}
	\end{definition}

	\begin{definition}\label{def:2.6}
		The point $x^* \in S$ is a stationary point of \eqref{qu:1} iff there exists a constant $\lambda>0$ such that
		\[
		x^*\in \arg\min_{u\in S} D_h(u, \nabla h^{*}\big(\nabla h(x^*) - \lambda \nabla f(x^*)\big)).
		\]
	\end{definition}

If the kernel function $h(x)=\frac{1}{2} \|x\|^2$, the stationary point conception recovers the $\lambda$-stationary point with respect to the classical hard-thresholding projection; see \cite[Definition 2.2]{zhou2022gradient}.
	

\section{Characterizations of Bregman Sparse Gradient Projection,  Local Minimizer and Stationary Point}	
	In this section, we present some characterizations of Bregman sparse gradient projection,  local minimizer and stationary point of problem \eqref{qu:1}, which will play a crucial theoretical tool to analyze the convergence of BGPNP.

The following theorem presents an equivalent characterization of local minimizer of \eqref{qu:1}.

	 \begin{theorem}\label{th2.1}
	 	If  $f$ is a convex function, then the vector $x^*$ is a local minimizer of \eqref{qu:1} if and only if it satisfies
	 \begin{equation*}
	 	\begin{split}
	 	\nabla_{\Gamma^{*}}f(x^*)=0,\,\,\text{if}\,\, \|x^*\|_0=s,\\
	 	\nabla f(x^*)=0,\,\,\text{if}\,\, \|x^*\|_0<s,
	 \end{split}
	\end{equation*}
	where $\Gamma^{*}=\operatorname{supp}(x^*)$.
	 \end{theorem}
	
	 \begin{proof}
	 	It follows from the proof of \cite[Theorem 10.1] {rockafellar2009variational}, \cite[Theorem 2.1] {zhou2022gradient} and \cite[Table 2]{raih2015solution}.
	 \end{proof}

	 The following theorem establishes the equivalence between the Bregman sparse gradient projection and a regularized first-order minimization problem.
	
	 \begin{theorem}\label{th:3.1}
	 	Let $h$ be a Legendre function, $f$ be a continuously differentiable function, $x\in \intr \dom h$ and $\lambda > 0$.  Define the Bregman gradient iteration
	 	\begin{align}\label{BG:1}
	 	x^+ = \nabla h^*(\nabla h(x) - \lambda \nabla f(x)).
	 	\end{align}
	 	Let the Bregman sparse gradient projection be
	 	\begin{align}\label{BSGP:1}
	 	x_s^+ \in \arg\min_{u\in S} D_h(u, x^+).
	 	\end{align}
	 	Then it holds that
	 	\[
	 	x_s^+ \in \arg\min_{u\in S}
	 	\left\{
	 	\langle \nabla f(x), u \rangle + \frac{1}{\lambda} D_h(u,x)
	 	\right\}.
	 	\]
	 	
	 \end{theorem}

	 \begin{proof}
	 		To derive the stated result, we proceed by proving
	 	that for any $s$-sparse vector $x_s$, i.e., $\|x_s\|_0 \le s$, the following inequality holds:
	 	\[
	 	\left\langle\nabla f\left(x\right), x^{+}_{s}\right\rangle+\frac{1}{\lambda} D_{h}\left(x^{+}_{s}, x\right) \leq \left\langle\nabla f\left(x\right), x_{s}\right\rangle+\frac{1}{\lambda} D_{h}\left(x_{s}, x\right).
	 	\]
	 	Then the above inequality is equivalent to the following inequality:
	 	\begin{equation}
	 		\begin{split}
	 			&\left\langle\nabla f\left(x\right), x^{+}_{s}-x^{+}\right\rangle+\frac{1}{\lambda}[ D_{h}\left(x^{+}_{s}, x\right)-D_{h}\left(x^{+}, x\right)] \\
	 			&\leq \left\langle\nabla f\left(x\right), x_{s}-x^{+}\right\rangle+\frac{1}{\lambda} [D_{h}\left(x_{s}, x\right)-D_{h}\left(x^{+}, x\right)]. \label{eq:thm1-2}
	 		\end{split}
	 	\end{equation}
	 	From the Bregman gradient iteration \eqref{BG:1}, we have
	 	\begin{equation*}
	 		x^{+}=\argmin_{u \in \mathbb{R}^{n}}\left\{\left\langle\nabla f\left(x\right), u\right\rangle+\frac{1}{\lambda} D_{h}\left(u, x\right)\right\}.
	 	\end{equation*}
	 	Applying the three-point identity of the Bregman distance, we get
	 	\begin{equation}
	 	\begin{split}
	 		&D_{h}\left(x^{+}_{s}, x\right)-D_{h}\left(x^{+}, x\right) \\
	 		&= D_{h}\left(x^{+}_{s}, x^{+}\right)+\left\langle\nabla h\left(x^{+}\right)-\nabla h\left(x\right), x^{+}_{s}-x^{+}\right\rangle,\label{eq:thm1-3}\\
	 \end{split}
 \end{equation}
 and
	 \begin{align}
	 		&D_{h}\left(x_{s}, x\right)-D_{h}\left(x^{+}, x\right) \nonumber\\
	 		&= D_{h}\left(x_{s}, x^{+}\right)+\left\langle\nabla h\left(x^{+}\right)-\nabla h\left(x\right), x_{s}-x^{+}\right\rangle. \label{eq:thm1-4}
	 	\end{align}
	 	Combining \eqref{eq:thm1-3} and \eqref{eq:thm1-4}, \eqref{eq:thm1-2} can be equivalently reformulated as
	 	\begin{equation}
	 		\begin{split}
	 			&\left\langle\nabla f\left(x\right), x^{+}_{s}-x_{s}\right\rangle+\frac{1}{\lambda}\big[ D_{h}\left(x^{+}_{s}, x^{+}\right)-D_{h}\left(x_{s}, x^{+}\right) \\
	 			&\quad +\left\langle\nabla h\left(x^{+}\right)-\nabla h\left(x\right), x^{+}_{s}-x^{+}\right\rangle \\
	 			&\quad -\left\langle\nabla h\left(x^{+}\right)-\nabla h\left(x\right), x_{s}-x^{+}\right\rangle\big] \leq 0. \label{eq:thm1-5}
	 		\end{split}
	 	\end{equation}
	 	Substituting $x^{+}=\nabla h^{*}(\nabla h(x)-\lambda \nabla f(x))$ into the left side of \eqref{eq:thm1-5} yields that
	 	\begin{align}
	 			&\left\langle\nabla f\left(x\right), x^{+}_{s}-x_{s}\right\rangle+\frac{1}{\lambda}\big[ D_{h}\left(x^{+}_{s}, x^{+}\right)-D_{h}\left(x_{s}, x^{+}\right) \nonumber\\
	 		&\quad +\left\langle\nabla h\left(x^{+}\right)-\nabla h\left(x\right), x^{+}_{s}-x^{+}\right\rangle-\left\langle\nabla h\left(x^{+}\right)-\nabla h\left(x\right), x_{s}-x^{+}\right\rangle\big]\nonumber\\
	 		&=\left\langle\nabla f\left(x\right), x^{+}_{s}-x_{s}\right\rangle+\frac{1}{\lambda}\big[ D_{h}\left(x^{+}_{s}, x^{+}\right)-D_{h}\left(x_{s}, x^{+}\right)-\lambda\left\langle \nabla f(x), x^{+}_{s}-x_{s}\right\rangle\big] \nonumber\\
	 		&= \frac{1}{\lambda}\big[ D_{h}\left(x^{+}_{s}, x^{+}\right)-D_{h}\left(x_{s}, x^{+}\right)\big]. \label{eq:thm1-6}
	 	\end{align}
	 	
	 	Since $x^{+}_{s}=\argmin_{u\in S} D_{h}\left(u, x^{+}\right)$, then for any $s$-sparse $x_{s}$,
	 	\[
	 	D_{h}\left(x^{+}_{s}, x^{+}\right) \leq D_{h}\left(x_{s}, x^{+}\right),
	 	\]
	 	i.e., $D_{h}\left(x^{+}_{s}, x^{+}\right)- D_{h}\left(x_{s}, x^{+}\right)\leq 0$,
	 	this together with \eqref{eq:thm1-6} implies that \eqref{eq:thm1-5} holds.
	 Therefore, the desired result is true by \eqref{eq:thm1-2} and \eqref{eq:thm1-6}, which completes the proof.
	 \end{proof}

	 Theorem \ref{th:3.1} provides the basic variational characterization of the Bregman projection iteration. This characterization is important because it allows us to compare the newly generated sparse point with any other feasible sparse vector.

For the convenience, we assume that the following Assumption \ref{assup:3.1} always  holds in the rest of this paper.

	\begin{assumption}\label{assup:3.1}
			Let $h$ be a Legendre function with $S \subseteq \intr \dom h$, $f$ be a continuously differentiable function and the pair of functions $(f,h)$ satisfy the LC condition with constant $L>0$.
		\end{assumption}

 We next present a three-point inequality, which is the main technical tool for proving descent and convergence of gradient type method in the sense of Bregman distance.
	
		\begin{lemma}
			Fix $x \in \mathrm{int}\mathrm{dom}\,h$ and  $\lambda > 0$. For any $u \in \mathrm{dom}\,h$, the following inequality holds:
		\begin{equation*}
			\left\langle\nabla f(x), x^{+}-u\right\rangle \leq \frac{1}{\lambda}\left(D_{h}(u, x)-D_{h}\left(u, x^{+}\right)-D_{h}\left(x^{+}, x\right)\right),
		\end{equation*}
where $x^{+}$ is defined by \eqref{BG:1}.
	\end{lemma}
		\begin{proof}
			It directly follows from the proof of \cite[Lemma 4.1, p.1081]{bauschke2019linear}. The proof is completed.
		\end{proof}

		\begin{lemma}\label{lem:3.2}
			For a given $x\in\operatorname{int}\operatorname{dom}h$
			and $\lambda>0$, let  $x^+$ be defined by \eqref{BG:1}
			and  $x_s^+$ be defined by \eqref{BSGP:1}.
			Then, for any $s$-sparse vector
			$u\in\operatorname{dom}h$, it holds that
			\begin{equation*}\label{eq:lemma32}
				\left\langle\nabla f(x),x_s^+-u\right\rangle
				\leq
				\frac{1}{\lambda}
				\left[
				D_h(u,x)-D_h(x_s^+,x)
				\right].
			\end{equation*}
		\end{lemma}
		
		\begin{proof}
			By the definition of the Bregman sparse projection,
			we have
			\begin{equation}\label{eq:projection_optimality}
				D_h(x_s^+,x^+)\leq D_h(u,x^+),
			\,\,\, \forall u\in S.
			\end{equation}
		Note that $	D_h(u,x^+)-D_h(x_s^+,x^+)
		=h(u)-h(x_s^+)
		-\left\langle
		\nabla h(x^+),u-x_s^+
		\right\rangle$.
			Together with $\nabla h(x^+)
			=\nabla h(x)-\lambda\nabla f(x)$,
			we obtain
			\begin{align*}
				&D_h(u,x^+)-D_h(x_s^+,x^+)\notag\\
			&=h(u)-h(x_s^+)
				-\left\langle
				\nabla h(x),u-x_s^+
				\right\rangle
				+\lambda\left\langle
				\nabla f(x),u-x_s^+
				\right\rangle.
			\end{align*}
			By the definition of the Bregman distance, one has
			\begin{align*}
				h(u)-h(x_s^+)
				-\left\langle
				\nabla h(x),u-x_s^+
				\right\rangle
				=D_h(u,x)-D_h(x_s^+,x).
			\end{align*}
			Consequently, we have
			\begin{align*}
				D_h(u,x^+)-D_h(x_s^+,x^+)
				=D_h(u,x)-D_h(x_s^+,x)+
				\lambda\left\langle
				\nabla f(x),u-x_s^+
				\right\rangle.
			\end{align*}
			Combining this identity with
			\eqref{eq:projection_optimality}, we obtain
			\begin{align*}
				\lambda\left\langle
				\nabla f(x),x_s^+-u
				\right\rangle
				\leq
				D_h(u,x)-D_h(x_s^+,x).
			\end{align*}
			Dividing both sides by $\lambda>0$ yields
			\[
			\left\langle
			\nabla f(x),x_s^+-u
			\right\rangle
			\leq
			\frac{1}{\lambda}
			\left(
			D_h(u,x)-D_h(x_s^+,x)
			\right).
			\]
			This completes the proof.
		\end{proof}

    The preceding lemma converts the optimality of the Bregman projection into an inequality involving the gradient of \(f\) and several Bregman distance terms. Combining this inequality with descent lemma yields the following fundamental inequality, which serves as the basement for the global convergence analysis.

    \begin{theorem}[Fundamental inequality]\label{thm:3.2}
    For a given $x\in\operatorname{int}\operatorname{dom}h$
	and $\lambda>0$,  let  $x^{+}_{s} \in \Pi_{S}(\nabla h^{*}(\nabla h(x) - \lambda \nabla f(x)))$. Then, for any $s$-sparse vector $u \in \mathrm{dom}h$, the following inequality holds:
\begin{align}\label{th:3.2:0}
	f\left(x^{+}_{s}\right) \leq  f(x)+\langle\nabla f(x), u-x\rangle+\frac{1}{\lambda} D_{h}(u, x) +\left(L-\frac{1}{\lambda}\right) D_{h}\left(x^{+}_{s}, x\right).
		\end{align}
	\end{theorem}
	
	\begin{proof}
	 By \cref{lem:3.2}, we have
		\begin{equation}
			\left\langle\nabla f(x), x^{+}_{s}-u\right\rangle \leq \frac{1}{\lambda}\left(D_{h}(u, x)-D_{h}\left(x^{+}_{s}, x\right)\right). \label{eq:thm2-1}
		\end{equation}
		According to \cref{lem:1}, one has
		\begin{equation}
			f(x^{+}_{s}) \leq f(x)+\langle\nabla f(x), x^{+}_{s}-x\rangle+LD_{h}(x^{+}_{s}, x). \label{eq:thm2-2}
		\end{equation}
		It thus follows from \eqref{eq:thm2-1} and \eqref{eq:thm2-2} that \eqref{th:3.2:0} holds.
		The proof is completed.
	\end{proof}

		\begin{lemma}
		\label{lem:3}
		Let $\lambda>0$ and $x\in \intr \dom h$. Define
		\begin{align}\label{lem:3.4:eq1}
			x^{+}_{\lambda,s}\in \Pi_S(\nabla h^{*}(\nabla h(x)-\lambda \nabla f(x))),\,\,
			x^{+}_{1,s}\in \Pi_S(\nabla h^{*}(\nabla h(x)-\nabla f(x))).
		\end{align}
		Then for any $\lambda \ge 1$, the following holds
		\[
		D_h(x^{+}_{\lambda,s}, x) \ge D_h(x^{+}_{1,s}, x).
		\]
		Moreover,  the uniform $\theta$-condition
		holds, i.e.,  there exists a function $\theta: \mathbb{R}_{++} \rightarrow  \mathbb{R}_{++}$  such that
		\[
		D_h(x^{+}_{\lambda,s}, x) \ge \theta(\lambda)\, D_h(x^{+}_{1,s}, x).
		\]
	\end{lemma}
	
	\begin{proof}
		It follows from \eqref{lem:3.4:eq1} that
		\[
		x^{+}_{\lambda,s}\in \argmin_{u\in S}\left\{\left\langle\nabla f\left(x\right), u\right\rangle+\frac{1}{\lambda} D_{h}\left(u, x\right)\right\}
		\]
		and
		\[
		x^{+}_{1,s}\in \argmin_{u\in S}\left\{\left\langle\nabla f\left(x\right), u\right\rangle+ D_{h}\left(u, x\right)\right\}.
		\]
		Then we can get the following inequalities
		\begin{align}
			\left\langle\nabla f\left(x\right), x^{+}_{\lambda,s}\right\rangle+\frac{1}{\lambda} D_{h}\left(x^{+}_{\lambda,s}, x\right) &\leq \left\langle\nabla f\left(x\right), x^{+}_{1,s}\right\rangle+\frac{1}{\lambda} D_{h}\left(x^{+}_{1,s}, x\right), \label{eq:lem3-1}\\
			\left\langle\nabla f\left(x\right), x^{+}_{1,s}\right\rangle+ D_{h}\left(x^{+}_{1,s}, x\right) &\leq \left\langle\nabla f\left(x\right), x^{+}_{\lambda,s}\right\rangle+ D_{h}\left(x^{+}_{\lambda,s}, x\right). \label{eq:lem3-2}
		\end{align}
		Combining \eqref{eq:lem3-1} with \eqref{eq:lem3-2}, we have
		\[
		\begin{aligned}
			\frac{1}{\lambda} [D_{h}\left(x^{+}_{\lambda,s}, x\right)-D_{h}\left(x^{+}_{1,s}, x\right)] &\leq \left\langle\nabla f\left(x\right), x^{+}_{1,s}-x^{+}_{\lambda,s}\right\rangle \\
			&\leq D_{h}\left(x^{+}_{\lambda,s}, x\right)-D_{h}\left(x^{+}_{1,s}, x\right).
		\end{aligned}
		\]
		Hence, when $\lambda\geq 1$,
		\[
		D_{h}\left(x^{+}_{\lambda,s}, x\right) \geq D_{h}\left(x^{+}_{1,s}, x\right).
		\]
		It also implies that the uniform $\theta$-condition holds. The proof is completed.
	\end{proof}

	The following result shows that every local minimizer of problem
		\eqref{qu:1} is also its strict local
			minimizer.
		
		\begin{theorem}
			Assume that $f$ is $l$-Legendre strongly convex
			with respect to the function $h$.
			If $x^*\in S$ is a local minimizer of problem
		\eqref{qu:1}, then $x^*$ is a strict local
			minimizer of problem \eqref{qu:1}.
		\end{theorem}
		
		\begin{proof}
			Let $x^*\in S$ be a local minimizer of
			problem \eqref{qu:1}, and denote
			$\Gamma^*=\operatorname{supp}(x^*)$.
			We split the proof into two cases.
			
			(C1) Let us consider the case $\|x^*\|_0=s$. Since $x^*$ is a local minimizer, then $	\nabla_{\Gamma^*}f(x^*)=0$. For any feasible point $x\in S$ sufficiently
			close to $x^*$, one has $	\operatorname{supp}(x)=\Gamma^*$ and so $\operatorname{supp}(x-x^*)  \subset\Gamma^*   $  Then
			\begin{align*}
				\langle\nabla f(x^*),x-x^*\rangle
				=\langle\nabla_{\Gamma^*}f(x^*),
				x_{\Gamma^*}-x_{\Gamma^*}^*\rangle
				=0.
			\end{align*}
By the $l$-Legendre strong convexity of $f$
with respect to $h$, we obtain
\begin{align*}
	f(x)
	&\geq f(x^*)
		+\langle\nabla f(x^*),x-x^*\rangle
				+lD_h(x,x^*)\\
				&=f(x^*)+lD_h(x,x^*).
			\end{align*}
Since $D_h(x,x^*)>0$ whenever $x\neq x^*$, we have
 \begin{align*}
 f(x)>f(x^*)
	\end{align*}
for all feasible point $x\neq x^*$ sufficiently
			close to $x^*$.

(C2) Let us consider the case $\|x^*\|_0<s$. Since $x^*$ is a local minimizer and
			$\|x^*\|_0<s$, then $\nabla f(x^*)=0$ and
			\begin{align*}
				f(x)
				&\geq f(x^*)
				+\langle\nabla f(x^*),x-x^*\rangle
				+lD_h(x,x^*)\\
				&=f(x^*)+lD_h(x,x^*).
			\end{align*}
			Since $D_h(x,x^*)>0$ for every $x\neq x^*$,
			we conclude that
			\[
			f(x)>f(x^*).
			\]
			Altogether, we conclude that $x^*$ is a strict local minimizer of
			problem \eqref{qu:1}. The proof is completed.
		\end{proof}

	The following result shows the relationship between stationary point and global minimizer  of problem \eqref{qu:1}.	

 \begin{theorem}\label{th:3.4}
 Assume that $f$ is $l$-Legendre strongly convex with respect to the function $h$. Let $x^*\in S$ be a stationary point of problem \eqref{qu:1} with  $\lambda>0$, i.e.,
 \[
 x^*\in\arg\min_{u\in S}
 D_h\left(
 u,\nabla h^*(\nabla h(x^*)-\lambda\nabla f(x^*))
 \right).
 \]
 If $\lambda>\frac1l$, then $x^*$ is the unique global minimizer of problem \eqref{qu:1}.
 \end{theorem}
\begin{proof}
Since $f$ is $l$-Legendre strongly convex with respect to $h$, then for every $u\in S$,
\[
	f(u)\ge f(x^*)
	+\langle\nabla f(x^*),u-x^*\rangle
	+lD_h(u,x^*).
\]
Since $x^*$ is a stationary point of problem \eqref{qu:1} associated with the step size $\lambda>0$,
\begin{align}\label{eq:3.4:1}
x^*\in\arg\min_{u\in S}
D_h\left(
u,\nabla h^*(\nabla h(x^*)-\lambda\nabla f(x^*))
\right).
\end{align}
By Theorem \ref{th:3.1}, this is equivalent to
\[
x^*\in\arg\min_{u\in S}
\left\{
\langle\nabla f(x^*),u\rangle
+\frac1\lambda D_h(u,x^*)
\right\}.
\]
Therefore, for every $u\in S$,
\[
\langle\nabla f(x^*),x^*\rangle
\le
\langle\nabla f(x^*),u\rangle
+\frac1\lambda D_h(u,x^*).
\]
Rearranging yields that
\begin{align}\label{eq:3.4:2}
	\langle\nabla f(x^*),u-x^*\rangle
	\ge-\frac1\lambda D_h(u,x^*).
\end{align}
Combining \eqref{eq:3.4:1} and \eqref{eq:3.4:2}, we obtain
\[
\begin{aligned}
	f(u)
	&\ge f(x^*)
	+\langle\nabla f(x^*),u-x^*\rangle
	+lD_h(u,x^*)\\
	&\ge f(x^*)+
	\left(l-\frac1\lambda\right)D_h(u,x^*).
\end{aligned}
\]
 Since $h$ is a Legendre function and $\lambda>\frac1l$, then $l-\frac1\lambda>0$ and
\[
D_h(u,x^*)>0,
\,\,\, \forall u\in S,\ u\ne x^*.
\]
Consequently, we have
\[
	f(u)>f(x^*),
	\,\,\, \forall u\in S,\ u\ne x^*.
\]
Therefore, $x^*$ is the unique global minimizer of problem \eqref{qu:1}.
This completes the proof.
\end{proof}

 It follows from the proof of Theorem \ref{th:3.4} that every stationary point of \eqref{qu:1} is its global minimizer when $f$ is a $l$-Legendre strongly convex function  with respect to $h$ and $\lambda \geq 1/l$.

\begin{theorem}\label{lem:3.3:0}
 If $x^*\in S$ is a stationary point of problem \eqref{qu:1} with constant $\bar{\lambda}>0$, then $x^*$ is also a stationary point of problem \eqref{qu:1} for all $\lambda \in (0,\bar{\lambda}]$.
\end{theorem}

 \begin{proof}
 	If $x^*\in S$ is a stationary point of problem \eqref{qu:1} with constant $\bar{\lambda}>0$, then
 		\[
 	x^*\in \arg\min_{u\in S} D_h(u, \nabla h^{*}\big(\nabla h(x^*) - \bar{\lambda} \nabla f(x^*)\big)).
 	\]
 	By Theorem \ref{th:3.1}, we have
 		\[
 x^* \in \arg\min_{u\in S}
 	\left\{
 	\langle \nabla f(x^*), u \rangle + \frac{1}{\bar{\lambda}} D_h(u,x^*)
 	\right\},
 	\]
 	which is equivalent to
 		\[
 	\langle \nabla f(x^*), u \rangle + \frac{1}{\bar{\lambda}} D_h(u,x^*)
 	\geq 	\langle \nabla f(x^*), x^* \rangle + \frac{1}{\bar{\lambda}} D_h(x^*,x^*), \,\,\, \forall \, u \in S.
 	\]
 	Then for each $u\in S$,
 	\[
 	\langle \nabla f(x^*), x^*-u \rangle \leq \frac{1}{\bar{\lambda}} D_h(u,x^*) \leq \frac{1}{\lambda} D_h(u,x^*), \,\,\,  \forall \, \lambda \in (0,\bar{\lambda}].
 	\]
 	Consequently, we conclude that for any $\lambda \in (0,\bar{\lambda}]$,
 		\[
 	x^* \in \arg\min_{u\in S}
 	\left\{
 	\langle \nabla f(x^*), u \rangle + \frac{1}{\lambda} D_h(u,x^*)
 	\right\},
 	\]
 	which is equivalent to
 		\[
 	x^*\in \arg\min_{u\in S} D_h(u, \nabla h^{*}\big(\nabla h(x^*) - \lambda \nabla f(x^*)\big)).
 	\]
 	Therefore the desired result holds. The proof is completed.
 \end{proof}

	\section{A Non-Euclidean Gradient Projection Newton Pursuit method}
	\label{sec:algorithm}

We now present a non-Euclidean Gradient Projection Newton Pursuit method for solving sparse constrained optimization problem \eqref{qu:1} via the Bregman distance. The method is also denoted as  BGPNP. The algorithm consists of two-stage. The first stage is a Bregman gradient projection step and the second stage is a Newton pursuit step restricted to the current active support. The fully convergence and locally linear convergence rate of BGPNP are established under some suitable conditions.

First, we denote $g^k=\nabla f(u^k)$,  $H_k=\nabla^2_{\Gamma_k\Gamma_k}f(u^k)$, where the index set $\Gamma_k$ is chosen and satisfies $\Gamma_k\supseteq \supp(u^k)$ and $|\Gamma_k|=s$.

	\begin{algorithm}[H]
		\caption{BGPNP: Bregman Gradient Projection Newton Pursuit}
		\label{alg:GPNP}
		\begin{algorithmic}[1]
			\Require Initial point $x^0$, parameters $\tau > 0$, $\sigma > 0$, $1 > \gamma > 0$, $\pi_0 >  \epsilon > 0$, $\kappa > 0$.
			\Ensure Solution $x^k$ of sparsity constrained optimization
			\State Set $k = 0$
			\While{$\pi_k > \epsilon$}
			\State \textbf{Bregman gradient projection}: Find the smallest integer $q_k = 0,1,\dots$ such that $\tau\gamma^{q_k}\leq\frac{1}{L+\sigma}$ and
			\[
			f(x^k(\tau\gamma^{q_k})) \leq f(x^k) - \sigma D_{h}(x^k(\tau\gamma^{q_k}),x^k),
			\]
			where \[
			x^k(\tau\gamma^{q_k}) \in \Pi_S(\nabla h^{*}(\nabla h(x^k)-\tau\gamma^{q_k} \nabla f(x^k))).
			\]
			\State Set $\alpha_k = \tau\gamma^{q_k}$, $u^k = x^k(\alpha_k)$, $x^{k+1}=u^{k}$.
			\State \textbf{Newton pursuit}: If $\supp(x^k) = \Gamma_k$ or $\norm{g^k} < \kappa$ (where $g^k = \nabla f(u^k), \Gamma_k \supseteq \supp(u^k), |\Gamma_k|=s$), compute $v^{k}$ by solving the system of equations:
			\[
				H_k(v^k_{\Gamma_k} - u^k_{\Gamma_k}) = -g^k_{\Gamma_k}, \quad v^k_{\overline{\Gamma}_k} = 0,
			\]
			where $H_k = \nabla^2_{\Gamma_k\Gamma_k}f(u^k)$.
			\State \quad If the solution $v^k$ satisfies
			\[
			f(v^k) \leq f(u^k) - \sigma D_{h}(v^k,u^k),
			\]
			\State \quad \quad set $x^{k+1} = v^k$.
			\State Compute $\pi_{k+1}$ (halting criterion) and set $k = k + 1$.
			\EndWhile
			\State Output $x^k$
		\end{algorithmic}
	\end{algorithm}
	
	In Algorithm \ref{alg:GPNP},
\[\pi_k =
\begin{cases}
	\|\nabla f(x^k)\|, & k < k_0, \\
	\max\big\{\operatorname{std}(f_k, f_{k-1}, \dots, f_{k-k_0}), \|\nabla f(x^k)\|\big\}, & k \ge k_0,
\end{cases}
\]
    where $f_k := f(x^k)$, $k_0$ is an integer (e.g., $k_0 = 5$ in the numerical experiments) and $\operatorname{std}(x)$ calculates the standard deviation of $x$. If $h(x)=\frac{\|x\|^2}{2}$, BGPNP reduces to the algorithm of \cite{zhou2022gradient}. In this sense, BGPNP provides an flexible framework for solving sparse constrained problem \eqref{qu:1} by choosing suitable kernel function $h$.

\subsection{Convergence Analysis}

	In this subsection, we establish the basic convergence properties, such as the bound of step size,  the non-increasing property for the values of objective function,  and asymptotic convergence of iteration points, of the BGPNP method under Assumption \ref{assup:3.1}. Moreover, the full convergence of iteration points is also established.

The following lemma shows that the Bregman gradient projection step guarantees a sufficient decrease of the objective function.

\begin{lemma}\label{lemm:3.5:001}
				For any $x\in S$, $\lambda>0$, and
				$
				x_s^+\in
				\Pi_S\left(
				\nabla h^*(\nabla h(x)-\lambda\nabla f(x))
				\right)$,
				it holds that
				\begin{equation*}\label{eq:descent35}
					f(x_s^+)
					\leq f(x)-
					\left(\frac{1}{\lambda}-L\right)
					D_h(x_s^+,x).
				\end{equation*}
				In particular, the sufficient decrease property
				holds whenever $0<\lambda L<1$.
				
\end{lemma}
		
\begin{proof}
			Taking $u=x\in S$ in \eqref{th:3.2:0}
		 of Theorem \ref{thm:3.2},
			we obtain
			\begin{align*}
				f(x_s^+)
			\leq f(x)-
				\left(\frac{1}{\lambda}-L\right)
				D_h(x_s^+,x).
			\end{align*}
			Since $0<\lambda L<1$, we have $\frac{1}{\lambda}-L>0$.
			Hence, the sufficient decrease property
			holds. The proof is completed.
      \end{proof}

		\begin{lemma}\label{lemm:3.5}
			Let $\{x^k\}$ be the sequence generated by the BGPNP.
			Then the following statements hold.
			\begin{enumerate}
				\item[\rm(i)]
				The step size $\alpha_k$ generated by the
				Bregman gradient projection step satisfies
				\[
				0<\underline{\alpha}
				\leq\alpha_k\leq\tau,
				\]
				where
				\begin{equation}\label{eq:lowerstepsize35}
					\underline{\alpha}
					=
					\min\left\{
					\tau,\frac{\gamma}{L+\sigma}
					\right\}>0.
				\end{equation}
				
				\item[\rm(ii)]
				The sequence $\{f(x^k)\}$ is non-increasing.
				
				\item[\rm(iii)]
				If $h$ is strongly convex with constant
				$\iota>0$, then
				\[
				\lim_{k\to\infty}\|u^k-x^k\|
				=
				\lim_{k\to\infty}\|x^{k+1}-x^k\|=0.
				\]
				
				\item[\rm(iv)]
				If $h$ is strongly convex with constant
				$\iota>0$, then every accumulation point
				of $\{x^k\}$ is a stationary point of
				problem \eqref{qu:1} with any  $\lambda \in (0,\alpha^*]$ for some $\alpha^*\in[\underline{\alpha},\tau]$.
			\end{enumerate}
		\end{lemma}
		
		\begin{proof}
			(i)
			By Algorithm \ref{alg:GPNP}, the step size $\alpha_k=\tau\gamma^{q_k}$, $0<\gamma<1$, where $q_k$ is the smallest nonnegative integer such that
		$\tau\gamma^{q_k}\leq\frac{1}{L+\sigma} $
	and
		\[
			f(x^k(\tau\gamma^{q_k}))
			\leq f(x^k)
			-\sigma D_h(x^k(\tau\gamma^{q_k}),x^k).
			\]
			Set $\bar{\alpha}=\frac{1}{L+\sigma}$.	For any $0<\lambda\leq\bar{\alpha}$, we have
			       $\frac{1}{\lambda}-L\geq\sigma$.
			Consequently, we conclude from Lemma \ref{lemm:3.5:001} that
			\[
			f(x^k(\lambda))
			\leq f(x^k)
			-\sigma D_h(x^k(\lambda),x^k),
			\]
where $x^k(\lambda)\in \Pi_{S}(\nabla h^{*}(\nabla h(x^k) - \lambda \nabla f(x^k)))$.
			Thus, the backtracking procedure is
			well defined.
			
			If $q_k=0$, then $\alpha_k=\tau\leq\bar{\alpha}$.
			If $q_k\geq 1$, by the minimality of $q_k$, we have 	$\tau\gamma^{q_k-1}> \bar{\alpha}$ or,
			\begin{align}\label{eq:30:0}
				f(x^k(\tau\gamma^{q_k-1})) > f(x^k) - \sigma D_{h}(x^k(\tau\gamma^{q_k-1}),x^k).
			\end{align}

			If $\tau\gamma^{q_k-1}> \bar{\alpha}$, multiplying both sides by $\gamma$, it follows that
			\[
			\tau\gamma^{q_k}
			>
			\gamma\bar{\alpha}
			=
			\frac{\gamma}{L+\sigma}.
			\]
			Therefore,
			\[
			\alpha_k
			>
			\underline{\alpha}
			=
			\min\left\{
			\tau,\,
			\frac{\gamma}{L+\sigma}
			\right\}
			>0.
			\]
			Since $\alpha_k=\tau\gamma^{q_k}\leq\tau$, we obtain
			\[
			0<\underline{\alpha}<\alpha_k\leq\tau,
			\,\,\, \forall k\geq0.
			\]
			
			If \eqref{eq:30:0} holds, then
			\[
			f(x^k(\tau\gamma^{q_k-1})) > f(x^k) - \sigma D_{h}(x^k(\tau\gamma^{q_k-1}),x^k).
			\]
			This together with Lemma \ref{lemm:3.5:001} yields that $\tau\gamma^{q_k-1}> \frac{1}{L}>  \frac{1}{L+\sigma}$. So, one has 		
			\[
			\alpha_k=\tau\gamma^{q_k}>\underline{\alpha}
			=
			\min\left\{
			\tau,\,
			\frac{\gamma}{L+\sigma}
			\right\}
			>0.
			\]
			Altogether, $\{\alpha_k\}$ is contained in the compact interval $[\underline{\alpha},\tau]$.

		(ii) By the sufficient decrease condition
			of Algorithm \ref{alg:GPNP}, we have
			\begin{equation}\label{eq:projdescent35}
				f(u^k)
				\leq f(x^k)
				-\sigma D_h(u^k,x^k).
			\end{equation}
			
			We distinguish two cases.
			
			\text{\bf Case 1}.
			If the Newton pursuit step is not
			accepted, then $x^{k+1}=u^k$.
			Therefore, one has
			\begin{align*}
				f(x^{k+1})
				=f(u^k)
				\leq f(x^k)
				-\sigma D_h(u^k,x^k)
				\leq f(x^k).
			\end{align*}

			\text{\bf Case 2}.
			If the Newton pursuit step is accepted,
			then $x^{k+1}=v^k$ and
			\[
			f(v^k)
			\leq f(u^k)
			-\sigma D_h(v^k,u^k).
			\]
			Combining this with
			\eqref{eq:projdescent35}, we obtain
			\begin{align*}
				f(x^{k+1})&=f(v^k)\leq f(u^k)
				-\sigma D_h(v^k,u^k)\nonumber\\
				&\leq f(x^k)
				-\sigma D_h(u^k,x^k)-\sigma D_h(v^k,u^k)\nonumber\\
				&\leq f(x^k).
			\end{align*}
			Altogether, we have
			\[
			f(x^{k+1})\leq f(x^k),
			\]
			and so, $\{f(x^k)\}$
			is non-increasing.
			
			(iii)
			 Since $f$ is bounded from below, (ii) shows that the non-increasing sequence  $\{f\left(x^{k}\right)\}$ converges and so,
			$$\sum_{k \geq 0} [f(x^k) - f(x^{k+1})]=f(x^0)-\lim_{k \to \infty} f(x^{k+1})<+\infty. $$
			Using the $\iota$-strong convexity of $h$, we obtain
			\[
			\begin{aligned}
				\text{\bf Case 1} \,\,&\Rightarrow\,\,
				2f(x^{k+1})\leq 2f(x^{k})-\sigma\iota\norm{u^{k}-x^{k}}^{2} ,\\
				\text{\bf Case 2} \,\,&\Rightarrow\,\,2f(x^{k+1})\leq 2f(x^{k})-\sigma\iota (\norm{u^{k}-x^{k}}^{2}+\norm{x^{k+1}-u^{k}}^{2})\\
				&\qquad\qquad\qquad\leq 2f(x^{k})-\frac{\sigma\iota}{2} \norm{x^{k+1}-x^{k}}^{2}.
			\end{aligned}
			\]
			Then
			\[
			\sum_{k\geq 0}\frac{\sigma\iota}{2}\norm{u^{k}-x^{k}}^{2}\leq f(x^0) - \lim_{k \to \infty} f(x^k) <+\infty ,\]
			which further leads to $\lim_{k \to \infty}\norm{u^{k}-x^{k}}=0$. Then $\lim_{k \to  \infty}\norm{x^{k+1}-x^{k}}=0$ holds in {\bf Case 1}.
			If {\bf Case 2} occurs, then
			\[
			\sum_{k\geq 0}\frac{\sigma\iota}{4}\norm{x^{k+1}-x^{k}}^{2} \leq f(x^0) - \lim_{k \to \infty} f(x^k) < +\infty,\]
			and so, $\lim_{k \to \infty}\norm{x^{k+1}-x^{k}}=0$.
			Consequently, we have
			\[
			\lim _{k \to \infty}\norm{u^{k}-x^{k}} = \lim _{k \to \infty}\norm{x^{k+1}-x^{k}} =0.
			\]

			(iv)
			Let $x^*$ be an arbitrary accumulation
			point of $\{x^k\}$. Without of loss generality, assume that $\lim _{k \to \infty}x^k=x^*$.  Using Algorithm \ref{alg:GPNP} yields that $x^*\in S$.
			By (iii), one has $\|u^{k}-x^{k}\|\to0$, and so, $	u^{k}\to x^*$.			
				It follows from (i) that for each $x^k$, there exists a corresponding step size $\alpha_k$ such that $\alpha_k\in[\underline{\alpha},\tau]$. Then there exist $\alpha^*>0$ and a convergent subsequence  $\{\alpha_{k_j}\}$ of $\{\alpha_k\}$ such that $\alpha_{k_j}\to\alpha^*\in[\underline{\alpha},\tau]$.
			For each $u \in S$, we have
			\[
			\lim_{k_j \to \infty} D_h(u, \nabla h^{*}\big(\nabla h(x^{k_j}) - \alpha_{k_j} \nabla f(x^{k_j})\big)) =
			D_h(u, \nabla h^{*}\big(\nabla h(x^*) - \alpha^* \nabla f(x^*)\big)).
			\]
For the brevity, set $x^+=\nabla h^{*}\left(\nabla h(x^*) - \alpha^* \nabla f(x^*)\right)$. Let us prove that $x^*$ is a stationary point.

			(c1) If $x^{k+1}$ is generated by Bregman projection step, we show the result by contradiction.
			Suppose that there exists $u\in S$ such that
			\[
			D_h(u, x^+) < D_h(x^*, x^+).
			\]
			By \cref{def:BD}, we have
			\begin{equation}
				\begin{split}\label{eq:30}
					&h(u)-h(x^*)-\langle\nabla h(x^*)-\alpha^*\nabla f(x^*), u-x^*\rangle \\
					&=D_h(u, x^*)- \alpha^*\langle\nabla f(x^*), x^*-u\rangle < 0,
				\end{split}
			\end{equation}
			because of $D_h(u, x^*)=h(u)-h(x^*)-\langle\nabla h(x^*), u-x^*\rangle$.
			By \cref{lem:3.2}, we have
			\[
			\langle \nabla f(x^k), x^{k+1} - u \rangle
			\leq \frac{1}{\alpha_k}
			\Big(
			D_h(u,x^k) - D_h(x^{k+1},x^k)
			\Big).
			\]
			Taking limit of both sides, we can derive that
			\[
			\alpha^* \langle\nabla f(x^*), x^*-u\rangle \leq D_h(u,x^*).
			\]
			This together with \eqref{eq:30} indicate that
			\[
			0\leq D_h(u, x^*) - \alpha^*\langle\nabla f(x^*), x^*-u\rangle < 0,
			\]
			which is a contradiction.
			Therefore, we obtain
			\begin{align}\label{eq:sp:0}
			x^*\in \arg\min_{u\in S} D_h(u, \nabla h^{*}\big(\nabla h(x^*) - \alpha^* \nabla f(x^*)\big)).
		\end{align}

			(c2) Let us consider the case where the Newton pursuit step is accepted, i.e., $x^{k+1}=v^k$. By (iv), we have
			\[
			\|u^k-x^k\|\to 0,
			\qquad
			\|x^{k+1}-x^k\|\to 0,
			\]
			which together with $x^k\to x^*$ yields that $u^k\to x^*$. By the triangle inequality, one has
			\[
			\|v^k-u^k\|
			\leq
			\|v^k-x^k\|
			+
			\|x^k-u^k\|.
			\]
			Since $v^k=x^{k+1}$, we obtain
			\[
			\|v^k-u^k\|
			\leq
			\|x^{k+1}-x^k\|
			+
			\|x^k-u^k\|\to 0.
			\]
		Then $\|v^k-u^k\|\to 0$.
			Due to $u^k\to x^*$, it follows that $v^k\to x^*$, which implies that $x^{k+1}\to x^*$.
			Consequently, whether the Newton pursuit step is accepted or not, the Bregman projection
			point $u^k$ and the actual next iterate $x^k$ have the same limit $x^*$.
			Therefore, the Newton pursuit step does not change the accumulation point of
			the sequence $\{x^k\}$. This combined with (c1) shows that $\{x^k\}$ converges to the stationary point $x^*$, satisfying  	\eqref{eq:sp:0}.
			It therefore follows from Theorem \ref{lem:3.3:0} that $x^*$ is the stationary point of problem \eqref{qu:1} with any $\lambda \in (0,\alpha^*]$.
			This completes the proof.
			\end{proof}

Now we present the convergence of the sequence
		$\{x^k\}$ generated by the BGPNP.

\begin{theorem}\label{thm:3.5}
	Assume that $h$ is strongly convex with constant
	$\iota>0$, and $f$ is Legendre strongly convex
	with respect to $h$ with constant $l>0$.
Let $\{x^k\}$ be an infinite sequence generated
	by the BGPNP.
	Then the following statements hold:
	\begin{enumerate}
		\item[\rm(i)] The sequence $\{x^k\}$ is bounded
		and admits at least one accumulation point.
		
		
		\item[\rm(ii)] The sequence
		$\{x^k\}$ fully converges to some stationary point of
		problem \eqref{qu:1}.
	\end{enumerate}
\end{theorem}

\begin{proof}
	(i) By Lemma \ref{lemm:3.5}(ii), the sequence $\{f(x^k)\}$ is non-increasing. Hence,
	\begin{equation}\label{eq:32:0}
		f(x^k)\leq f(x^0),
		\,\,\, \forall k\geq 0.
	\end{equation}
	Since $f$ is $l$-Legendre strongly convex
	with respect to the $\iota$-strongly convex
	function $h$, for any $x,y\in S$,
	we have
	\begin{align}\label{eq:33-1}
		f(x)
		&\geq f(y)
		+\langle\nabla f(y),x-y\rangle
		+lD_h(x,y)
		\nonumber\\
		&\geq f(y)
		+\langle\nabla f(y),x-y\rangle
		+\frac{l\iota}{2}\|x-y\|^2.
	\end{align}
	Due to $l\iota>0$, the quadratic term dominates
	the linear term as $\|x\|\to\infty$.
	Therefore, \eqref{eq:33-1} implies that $f$ is coercive.
	Define the level set
	\[
	\mathcal{L}_0
	:=
	\{x\in S:f(x)\leq f(x^0)\}.
	\]
	Since $f$ is continuous and coercive,
	and $S$ is closed, the level set
	$\mathcal{L}_0$ is compact.	Moreover, \eqref{eq:32:0} implies that $x^k\in\mathcal{L}_0$, for all $k\geq0$.
	Consequently, $\{x^k\}$ is bounded and
	admits at least one accumulation point.

	(ii) Let $\Omega$ denote its set of accumulation points of the sequence $\{x^k\}$. Then $S\supset \Omega\neq \emptyset$ by the boundedness of the sequence $\{x^k\}$. For any  $x^* \in \Omega$, it follows from Lemma \ref{lemm:3.5}(iv) that $x^*$ is a stationary point of problem \eqref{qu:1}.
 Hence, there exists $\lambda>0$ such that
   	\[
   	x^*\in\arg\min_{u\in S}
   	D_h\left(
   	u,\nabla h^*(\nabla h(x^*)-\lambda\nabla f(x^*))
   	\right).
   	\]
   	By Theorem \ref{th:3.1}, this is equivalent to
   	\[
   	x^*\in\arg\min_{u\in S}
   	\left\{
   	\langle\nabla f(x^*),u-x^*\rangle
   	+\frac{1}{\lambda}D_h(u,x^*)
   	\right\}.
   	\]
   	
  Set $\Gamma^*= \operatorname{supp}(x^*)$ and
   	 $	\Phi(u)=
   	\langle\nabla f(x^*),u-x^*\rangle
   	+\frac{1}{\lambda}D_h(u,x^*)$.	
   	Then
   	\begin{align} \label{eq:th3.5:1}
   		  	\Phi(u)\geq\Phi(x^*)=0,
   	\,\,\, \forall u\in S.
   \end{align}
   	
   	We distinguish the following two cases.
   	
   	\text{(c1)} Suppose that $\|x^*\|_0=s$.
   	For any $i\in\Gamma^*$, let $u(t)=x^*+te_i$, where $e_i$ is a unit vector of $\mathbb{R}^n$ and its $i$-th component is $1$.
   	Then $\operatorname{supp}(u(t))\subseteq\Gamma^*$
   	and $u(t)\in S$. Set $\psi_i(t)=\Phi(u(t))$. By \eqref{eq:th3.5:1}, we obtain  $$\psi_i(t)=\Phi(u(t))\geq\Phi(x^*)=\psi_i(0)=0, \,\,\, \forall t \in\mathbb{R},$$
   which implies that $t=0$ is a minimizer of $\psi_i$, and so $\psi_i'(0)=0$.
   	Note that
   	\[
   	\begin{aligned}
   		\psi_i'(t)
   		=\langle\nabla f(x^*),e_i\rangle
   		+\frac{1}{\lambda}
   		\left\langle
   		\nabla h(u(t))-\nabla h(x^*),e_i
   		\right\rangle.
   	\end{aligned}
   	\]
   	Taking $t=0$, we obtain
   	\[
   	0=\psi_i'(0)
   	=\langle\nabla f(x^*),e_i\rangle
   	=\frac{\partial f(x^*)}{\partial x_i}.
   	\]
   By the arbitrariness of $i\in\Gamma^*$, it follows that $\nabla_{\Gamma^*}f(x^*)=0$.

   	\text{(c2)} Suppose that $\|x^*\|_0<s$.
   	For any $i\in\{1,\ldots,n\}$, let $	u(t)=x^*+te_i$.
   	Then $\|u(t)\|_0
   	\leq\|x^*\|_0+1
   	\leq s$. So, $u(t)\in S$ and $$\Phi(u(t))\geq\Phi(x^*)=0,\,\,\, \forall t\in\mathbb{R}.$$
   	Set $\psi_i(t)=\Phi(u(t))$. So, $t=0$ is a minimizer of
   	$\psi_i(t)$, which implies
   	\[
   	0=\psi_i'(0)
   	=\frac{\partial f(x^*)}{\partial x_i}.
   	\]
   	By the arbitrariness of $i\in\{1,\ldots,n\}$, we conclude that

   $\nabla f(x^*)=0$.
   	
   	 Moreover, one has $\nabla_{\Gamma^*}f(x^*)=0$.
   From the arbitrariness of $x^*$, we conclude that $\nabla_{\operatorname{supp}(x^*)}f(x^*)=0$ for all $x^*\in \Omega$.

Since \(f\) is \(l\)-Legendre strongly convex with respect to the \(\iota\)-strongly convex function \(h\), we have
\begin{align*}
   		f(x)\geq  f(y)
   		+\langle\nabla f(y),x-y\rangle
   		+\frac{l\iota}{2}\|x-y\|^2,
   \end{align*}
which implies that $f$ is strongly convex and coercive. Then there exists at most one stationary point on each coordinate subspace $\mathbb{R}^{n}_{\bar{\Gamma}}=\{x\in \mathbb{R}^{n}\,:\, x_{\bar{\Gamma}}=0\}$, where $|\Gamma|\leq s$.
 Since there are only finitely many index sets $\Gamma$ satisfying $|\Gamma|\le s$, then
 the set of stationary points of problem \eqref{qu:1} is finite.
From Lemma \ref{lemm:3.5}(iv), we conclude that every accumulation point of $\{x^k\}$ is stationary point and so, $|\Omega|$ is finite.

Let us show that $\Omega$ is a singleton.  Suppose that there exist $x^*, y^*\in\Omega $ such that $x^*\neq y^*$.
 Since $\Omega$ is finite, then
\begin{equation*}\label{dist:AP}
\min_{y\in\Omega\setminus\{x^*\}}
\|y-x^*\|>0.
\end{equation*}
Set $d=\min_{y\in\Omega\setminus\{x^*\}}
\|y-x^*\| $ and
 $\delta=\frac{ d}{4}$.
  Then $\|y^*-x^*\|\ge d=4\delta$ and
\begin{align}\label{sing}
B(x^*,2\delta)\cap\Omega=\{x^*\},
\end{align}
where $B(x^*,2\delta)$ is a closed ball centred at $x^*$ with radius $2\delta$.

Since $x^*\in \Omega $ is an accumulation point,  there exists a subsequence $\{x^p\}$ of $\{x^k\}$ such that $x^p\rightarrow x^*$. Then
there exists a positive integer $K_{1}$ such that
\begin{align*}
\|x^p-x^*\|<\delta,\,\, \forall\, p\ge K_{1} .
\end{align*}
Similarly, there exists a subsequence $\{x^q\}$ of $\{x^k\}$ such that $x^q\rightarrow y^*$.
For every sufficiently large $p$, there exists  $q>p$ such that $ \|x^q- y^*\|\leq \delta $ and
\begin{align}\label{eq:13:01}
\|x^q-x^*\|\geq  \|y^* -x^* \|- \|x^q- y^*\|  >2\delta.
\end{align}
In turn, Lemma \ref{lemm:3.5}(iii) yields that
\begin{align*}
\lim_{k\to\infty}\|x^{k+1}-x^k\|=0,
\end{align*}
which implies that there exists a positive integer $K$ such that
\begin{align}\label{eq:14:1}
\|x^{k+1}-x^k\|<\delta,
\,\,\forall\, k\ge K.
\end{align}
Choose $p\ge \max \{K, K_{1} \}$, and let $q>p$ be the smallest index such that \eqref{eq:13:01} holds.
So, one has
\begin{align}\label{eq:15:1}
\|x^{q-1}-x^*\|\le2\delta.
\end{align}
Moreover, we have
\begin{align}\label{eq:16:1}
	\|x^{q-1}-x^*\|
	\ge\|x^q-x^*\|-\|x^q-x^{q-1}\|
	>2\delta-\delta
	=\delta,
\end{align}
where the last inequality results from \eqref{eq:13:01} and \eqref{eq:14:1}.
It follows from both \eqref{eq:15:1} and \eqref{eq:16:1} that
\begin{align}\label{eq:17:1}
	\delta<\|x^{q-1}-x^*\|\le2\delta.
\end{align}
In particular,
 we can choose a subsequence $\{x^{p_{j}}\}$ of $\{x^p\}$ satisfying the strictly increasing sequence of indices $\{p_j\}$ with $p_j\to\infty$.
 Then,
for each sufficiently large $p_j$,  $\|x^{p_j}-x^*\|<\delta$ and
there exists a subsequence $\{x^{q_{j}}\}$ of $\{x^q\}$ satisfying that $q_j$ is the corresponding smallest index such that $q_j>p_j$  and  $\|x^{q_{j}}-x^*\|>2\delta$.  Then
  $q_j-1\ge p_{j}\to\infty $
and $x^{q_{j}-1}$ satisfies \eqref{eq:17:1}. It infers that
  $x^{q_{j}-1}\in \mathcal{A}$,
where
   $\mathcal{A} = \left\{
x\in\mathbb{ R}^n \,:\,
\delta\le\|x-x^*\|\le2\delta \right\}$.
 Taking into account that $x^q\rightarrow y^*$,
 we deduce that $x^{q_{j}-1} \rightarrow y^*$. Then $y^*\in \mathcal{A}$
 because of the compactness of $\mathcal{A}$.
Moreover, $y^*\in \mathcal{A} \subset B(x^*,2\delta)$.
Consequently, $y^*\in B(x^*,2\delta) \cap\Omega$,
which contradicts \eqref{sing}.
Therefore, $\Omega$ is a singleton. It then implies that $\{x^k\}$ admits a unique accumulation point $x^{\ast}$, which is a stationary point of problem \eqref{qu:1}. Namely, the sequence
	$\{x^k\}$ fully converges to $x^{\ast}$.
	This completes the proof.
\end{proof}

\subsection{Linear Convergence}

    In this subsection, we analyze the linear convergence of the Bregman gradient projection without Newton pursuit step.

We firstly present the linear convergence in the sense of objective function for the Bregman gradient projection without Newton pursuit step.

\begin{theorem}[Linear convergence rate]\label{thm:3.3}
Let $\{x^k\}$ be the sequence generated
	by BGPNP and $x^{\ast}$ be an accumulation point of $\{x^k\}$ with $\|x^{\ast}\|_{0}=s$.
Assume that all conditions of Theorem \ref{thm:3.5} hold and $\nabla h$ is Lipschitz continuous
with constant $\varsigma>0$ on a convex neighborhood containing
$x^*$ of $\mathbb{R}^{n}_{\Gamma^{\ast}}$.
Then there exists an integer $K\geq 0$ such that
\begin{align}\label{linear:Obj}
f(x^k)-f(x^*)
\leq
(1-c)^{k-K}
\big(f(x^K)-f(x^*)\big),
\,\,\forall\, k\geq K,
\end{align}
i.e., the  sequence of objective function values converges
linearly to $f(x^*)$, where
\begin{align}\label{Parameters}
\Gamma^* =\operatorname{supp}(x^*),\,\,
 c=
\min\left\{
\frac12,
\frac{\sigma l\iota^2\underline{\alpha}^{\,2}}{M^2}
\right\}\,\,\,
\mbox{and}\,\,
\underline{\alpha}
=
\min\left\{
\tau,\frac{\gamma}{L+\sigma}
\right\}.
\end{align}
\end{theorem}

\begin{proof}
 It follows from Theorem \ref{thm:3.5}   that $\{x^k\}$ fully converges to $x^*$,
  which is a stationary point of
		problem \eqref{qu:1}.  Again from the proof of Theorem \ref{thm:3.5}(ii), we obtain
       $\nabla_{\Gamma^*}f(x^*)=0$.
It follows from Lemma \ref{lemm:3.5}(iii) that $u^k\to x^*$. Then $\Gamma^* \subset \operatorname{supp}(u^k) \subset \Gamma_k$ and $\Gamma^* \subset \operatorname{supp}(x^k) $ for all sufficiently large $k$.
 Due to $\|x^*\|_0=s$ and $x^k,u^k\in S$, for sufficiently large $k$, we have
   	\begin{align*}
   	\operatorname{supp}(x^k)=\operatorname{supp}(u^k)=\Gamma_k=\Gamma^*.
   	\end{align*}
So,  $x^k,u^k\in \mathbb{R}^{n}_{\Gamma^{\ast}}$  for all sufficiently large $k$.
In turn, we deduce from  Algorithm \ref{alg:GPNP} and
Theorem \ref{th:3.1}   that
\begin{align*}
u^k\in\arg\min_{u\in S}
\left\{
\langle\nabla f(x^k),u\rangle
+\frac{1}{\alpha_k}D_h(u,x^k)
\right\}.
\end{align*}
In view of
$\mathbb{R}^{n}_{\Gamma^{\ast}}\subseteq S$, we have
\begin{align}\label{eq:sub:optim:0}
u^k\in\arg\min_{u\in \mathbb{R}^{n}_{\Gamma^{\ast}}}
\left\{
\langle\nabla f(x^k),u\rangle
+\frac{1}{\alpha_k}D_h(u,x^k)
\right\},
\end{align}
 for all sufficiently large $k$. It follows from the convexity of $f$ and $h$ that for all sufficiently large $k$, \eqref{eq:sub:optim:0} can be equivalently characterized by
\begin{align*}
\nabla_{\Gamma^*}f(x^k)
+\frac{1}{\alpha_k}
\left(
\nabla_{\Gamma^*}h(u^k)
-\nabla_{\Gamma^*}h(x^k)
\right)=0,
\end{align*}
and so,
 $\alpha_k\nabla_{\Gamma^*}f(x^k)= \nabla_{\Gamma^*}h(x^k)-\nabla_{\Gamma^*}h(u^k)$.
Since $\nabla h$ is Lipschitz continuous
with constant $\varsigma>0$ on a convex neighborhood containing
$x^*$ of $\mathbb{R}^{n}_{\Gamma^{\ast}}$, then for any sufficiently large $k$,
\begin{align*}
\|u^k-x^k\|
 \geq
 \frac{\alpha_k}{\varsigma}
   \|\nabla_{\Gamma^*}f(x^k)\|.
\end{align*}
By the strong convexity of $h$ with constant
$\iota>0$, we obtain
\begin{align*}
D_h(u^k,x^k)
\geq
\frac{\iota}{2}\|u^k-x^k\|^2,
\end{align*}
Then
 $ D_h(u^k,x^k) \geq
\frac{\iota\alpha_k^2}{2\varsigma^2}
\|\nabla_{\Gamma^*}f(x^k)\|^2 $.
Using Lemma \ref{lemm:3.5}(i) yields that
\begin{equation}\label{eq:local_bregman_lower}
D_h(u^k,x^k)
\geq
\frac{\iota\underline{\alpha}^{\,2}}{2\varsigma^2}
\|\nabla_{\Gamma^*}f(x^k)\|^2.
\end{equation}

Since $f$ is Legendre strongly convex with
respect to $h$ with constant $l>0$, then  $f$ is $l\iota$-strongly
convex. Then $\min_{x\in \mathbb{R}^{n}_{\Gamma^*}}f(x)$ admits a unique solution.
This together with $\nabla_{\Gamma^*}f(x^*)=0$ show that
\begin{align}\label{eq:Optim:1}
x^* =
\arg\min_{x\in \mathbb{R}^{n}_{\Gamma^*}}f(x).
\end{align}
It follows from the strong convexity of  $f$ that the Polyak-Lojasiewicz inequality holds:
\begin{align*}
f(x)-f(x^*)
\leq
\frac{1}{2l\iota}
\|\nabla_{\Gamma^*}f(x)\|^2,\,\, \forall\, x\in \mathbb{R}^{n}_{\Gamma^*} .
\end{align*}
So, one has
\[
\|\nabla_{\Gamma^*}f(x^k)\|^2
\geq
2l\iota\big(f(x^k)-f(x^*)\big)
\]
for all sufficiently large $k$.
Moreover,
\eqref{eq:local_bregman_lower} implies that
\begin{align}\label{eq:local_error_bound}
D_h(u^k,x^k)
\geq
\frac{\iota\underline{\alpha}^{\,2}}{2\varsigma^2}
\|\nabla_{\Gamma^*}f(x^k)\|^2
\geq
\frac{l\iota^{2}\underline{\alpha}^{2}}{\varsigma^2}
\big(f(x^k)- f(x^*)\big).
\end{align}

Note that
 $f(x^{k+1}) \leq
f(x^k)-\sigma D_h(u^k,x^k)$.
Then
\begin{align*}
f(x^{k+1})-f(x^*)
&\leq f(x^k)-f(x^*)
-\sigma D_h(u^k,x^k)\\
&\leq
\left(
1- \frac{\sigma  l\iota^{2}\underline{\alpha}^{2}}{\varsigma^2}
\right)
\big(f(x^k)-f(x^*)\big).
\end{align*}
where the second inequality results from \eqref{eq:local_error_bound}. Set
 $c=
\min\left\{
\frac{1}{2},
\frac{\sigma  l\iota^{2}\underline{\alpha}^{2}}{\varsigma^2}
\right\}$.
Then $0<c<1$, and thus
\begin{align*}
0\leq f(x^{k+1})-f(x^*)
\leq
(1-c)\big(f(x^k)-f(x^*)\big)
\end{align*}
for all sufficiently large $k$. Therefore there exists a sufficiently large
 $K$ such that \eqref{linear:Obj} holds.
The proof is completed.
\end{proof}

The following result presents the linear convergence of the sequence $\{x^k\}$.
	
	\begin{theorem}[Linear convergence rate]\label{thm:3.7}
		Let $\{x^k\}$ be the sequence generated
		by the BGPNP. Assume that all conditions of Theorem \ref{thm:3.3} hold. Then there exists an integer $K\geq 0$ such that
         \begin{align*}
		\|x^{k+1}-x^k\|
		\leq
		(\sqrt{1-c})^{k-K} \sqrt{\frac{2(f(x^{K}) - f(x^*))}{\sigma \iota}},\,\,\forall\, k\geq K,
		\end{align*}
			where  $c\in (0,1)$ is defined by \eqref{Parameters}.
		\end{theorem}

\begin{proof}
	From the proof of Theorem \ref{thm:3.3}, we obtain that $x^k,u^k\in \mathbb{R}^{n}_{\Gamma^{\ast}}$  for all sufficiently large $k$.
 Then there exists a integer $K\geq 0$ such that
   \begin{align*}
f(x^k) \geq
\min_{x\in \mathbb{R}^{n}_{\Gamma^*}}f(x)=f(x^*),\,\,\forall\,k\geq K,
\end{align*}
where the inequality holds by \eqref{eq:Optim:1}.
It follows from Algorithm \ref{alg:GPNP} that
	\begin{align*}
		f(u^k)
	\leq f(x^k)-\sigma D_h(u^k,x^k).
	\end{align*}
This together with $u^k=x^{k+1}$ yields that
		\begin{align}\label{eq:3.7:1}
\sigma D_h(x^{k+1},x^k)
	\leq f(x^k) - f(x^{k+1}) \leq  f(x^k) - f(x^*),\,\,\forall\,k\geq K.
		\end{align}
	By the strong convexity of $h$, we have
	\[
	\frac{ \iota}{2} \|x^{k+1}-x^k\|^2
	\leq  D_h(x^{k+1},x^k).
	\]
	Consequently, \eqref{eq:3.7:1} implies that
	\begin{align}\label{eq:3.7:2}
	\frac{ \iota \sigma}{2} \|x^{k+1}-x^k\|^2
	& \leq f(x^k) -f(x^*) \nonumber\\
	& \leq (1-c)^{k-K}\big(f(x^K)- f(x^*)\big),\,\,\forall\, k\geq K,
		\end{align}
		where the last inequality follows from
	Theorem \ref{thm:3.3}.
	Therefore the desired result holds by \eqref{eq:3.7:2}. The proof is completed.
	\end{proof}

	\begin{remark}
	If $h(x)=\frac{\|x\|^2}{2}$, BGPNP reduces to the GPNP \cite{zhou2022gradient}. However, the linear convergence rate for the GPNP \cite{zhou2022gradient} is not established. So, Theorems \ref{thm:3.3} and \ref{thm:3.7} are original even if $h(x)=\frac{\|x\|^2}{2}$. Besides, if $s=n$, the problem \eqref{qu:1} reduces to an unconstrained optimization problem. In this case, the linear convergence results do not involve the gradient dominated condition and $\theta$ uniform condition, which are required in \cite[Theorem 4.1 and Proposition 4.1]{bauschke2019linear}. As mentioned in \cite{bauschke2019linear}, is not easy to verify the $\theta$ uniform condition, which also confined the application of \cite[Theorem 4.1 and Proposition 4.1]{bauschke2019linear}.
	\end{remark}

\section{Locally Quadratic Convergence}
		
	In this section, we investigate the local quadratic convergence behavior of BGPNP under suitable conditions.

	The following theorem shows that, under suitable conditions, the Newton pursuit step is eventually accepted and the sequence converges quadratically to the solution.

   \begin{theorem}[Locally quadratic convergence]\label{th:4.1}
   	Assume that $h$ is strongly convex with
   	constant $\iota>0$, and $f$ is Legendre strongly convex
   	with respect to $h$ with constant $l>0$. Let $\{x^k\}$ be an infinite sequence generated by the BGPNP, and let $x^k$ converge to $x^*$. Assume $f$ is locally $s$-RHLC at $x^*$ with constant $M_s^*>0$, and let
   	$0<\sigma<l$. Then the following statements hold.
   	
   	\begin{enumerate}
   		\item[\rm(i)]
   		For all sufficiently large $k$, the support
   		of $x^*$ is identified in the sense that
   		\begin{equation*}\label{eq:support_identification}
   			\begin{cases}
   				\operatorname{supp}(x^*)
   				\subseteq
   				\operatorname{supp}(x^k)\cap\Gamma_k,
   				& \|x^*\|_0<s,\\[2mm]
   				\operatorname{supp}(x^*)
   				=
   				\operatorname{supp}(x^k)
   				=
   				\Gamma_k,
   				& \|x^*\|_0=s.
   			\end{cases}
   		\end{equation*}
   		
   		\item[\rm(ii)]
   		The Newton pursuit step is activated and
   		accepted for all sufficiently large $k$.
   		
   		\item[\rm(iii)]
   		There exist constants $C>0$ and $K>0$
   		such that
   		\begin{equation*}\label{eq:quadratic_convergence}
   			\|x^{k+1}-x^*\|
   			\leq C\|x^k-x^*\|^2,
   			\,\,\,\forall k\geq K,
   		\end{equation*}
   		namely, $\{x^k\}$ converges
   		Q-quadratically to $x^*$.
   	\end{enumerate}
   \end{theorem}

   \begin{proof}
   It follows from Lemma \ref{lemm:3.5}(iii) and  $x^k\to x^*$ that $u^k\to x^*$. Let $\Gamma^*=\operatorname{supp}(x^*)$.
   	
   	(i) If $\|x^*\|_0=s$, since $x^k,u^k\in S$, then for sufficiently large $k$, we have
   	\[
   	\operatorname{supp}(x^k)=\operatorname{supp}(u^k)=\Gamma^*.
   	\]
   	Since $
   	\Gamma_k\supseteq\operatorname{supp}(u^k), |\Gamma_k|=s$, then $\Gamma_k=\Gamma^*$ and so, $	\operatorname{supp}(x^*)=\operatorname{supp}(x^k)=
	\Gamma_k.$

	If $\|x^*\|_0<s$, $x^k\to x^*$ and $u^k\to x^*$ implies that for sufficiently large $k$, $\Gamma^*
   	\subseteq\operatorname{supp}(x^k)$ and $\Gamma^*
   	\subseteq\operatorname{supp}(u^k)
   	\subseteq\Gamma_k$. Consequently, one has
   	\[
   	\Gamma^*
   	\subseteq
   	\operatorname{supp}(x^k)\cap\Gamma_k.
   	\]

   	(ii) We firstly show that the the Newton pursuit step
   	is well defined and  eventually activated.
   		By Lemma \ref{lemm:3.5}(iv), $x^*$ is a stationary point of problem \eqref{qu:1}.
%
   By the proof of Theorem \ref{thm:3.5} (ii), we deduce that
\begin{equation*}\label{eq:support_identification}
   			\begin{cases}
   					\nabla f(x^*)=0,
   				& \|x^*\|_0<s,\\[1mm]
   				\nabla_{\Gamma^*}f(x^*)=0,
   				& \|x^*\|_0=s.
   			\end{cases}
   		\end{equation*}
Together with (i) yields
   	that there exists an integer $K>0$ such that
   	\begin{equation}\label{eq:restricted_stationarity}
   		\nabla_{\Gamma_k}f(x^*)=0,
   		\,\,\, \forall k\geq K.
   	\end{equation}

   	By the $l$-Legendre strong convexity
   	of $f$ with respect to the
   	$\iota$-strongly convex function $h$, we have
   	\begin{align*}
   		f(x) \geq f(y)
   		+\langle\nabla f(y),x-y\rangle
   		+\frac{l\iota}{2}\|x-y\|^2,\,\,\forall \, x,y\in \mathbb{R}^{n},
   \end{align*}
    which implies that $\nabla^2 f(y)\succeq l\iota \mathbb{I}$ because of $f\in C^2$, where $\mathbb{I}$ is identity matrix. Therefore, for any index set
   	$\Gamma$ with $|\Gamma|\leq s$, one has
   	\[
   	\nabla^2_{\Gamma\Gamma}f(y)
   	\succeq l\iota \mathbb{I}.
   	\]
   	According to $H_k=\nabla^2_{\Gamma_k\Gamma_k}f(u^k)$, we obtain
   	\begin{equation}\label{eq:Hessian_inverse_41}
   		\lambda_{\min}(H_k)\geq l\iota>0,
   		\,\,\,
   		\|(H_k)^{-1}\|\leq\frac{1}{l\iota}.
   	\end{equation}
   	Hence, the the Newton pursuit step
   	is well defined.

   	Let us now show that the Newton pursuit step
   	is eventually activated. If $\|x^*\|_0<s$, then $\nabla f(x^*)=0$.
   	 Since $\nabla f$ is continuous and $u^k\to x^*$, then for all sufficiently large $k$, such that $\|\nabla f(u^k)\|<\kappa$, where $\kappa>0$ is the same as in Algorithm \ref{qu:1}. If $\|x^*\|_0=s$, (i) implies that for all sufficiently large $k$, $
   	 \operatorname{supp}(x^k)=\Gamma_k=\Gamma^*$.
   	 Consequently, the Newton pursuit step
   	is eventually activated.
   	
   	We next prove that the Newton pursuit step is
   		accepted for all sufficiently large $k$.
   	From Algorithm \ref{qu:1}, we have
   	\begin{equation}\label{hk:0}
   		H_k(v^k-u^k)_{\Gamma_k}
   		=
   		-\nabla_{\Gamma_k}f(u^k),\,\,\, v^k_{\overline{\Gamma}_k} = 0.
   	\end{equation}
   	Let $d^k=v^k-u^k$. Using \eqref{eq:restricted_stationarity} shows that
   \begin{equation}\label{eq:dk1}
   	\begin{aligned}
   		\|d^k\|	=\left\|
   		H_k^{-1}\nabla_{\Gamma_k}f(u^k)\right\|
   		=
   		\left\|	H_k^{-1}\bigl(\nabla_{\Gamma_k}f(u^k)-\nabla_{\Gamma_k}f(x^*)\bigr)\right\|.
   	\end{aligned}
   	   \end{equation}
   	Since $f\in C^2$, its Hessian $\nabla^2 f(\cdot)$ is
   	locally bounded around $x^*$.
   	Hence, for all $y$ in a sufficiently small
   	neighborhood of $x^*$, there exists $M_f>0$ such that
   	\begin{align}\label{eq:Mf}
   	\left\|\nabla^2_{\Gamma_k, \Gamma_k}f(y)\right\| \leq \|\nabla^2 f(y)\|\leq M_f.
   	\end{align}
 For any sufficiently large $k$, set $x_t^k = t x^k + (1-t)x^* \in S$ and $u_t^k = t u^k + (1-t)x^* \in S$,  $t \in [0,1]$. By the Mean Value Theorem, one has
   	\begin{align}\label{eq:so-3:00}
   	 \|\nabla f(x^k) - \nabla f(x^*)\|
    &= \left\| \int_0^1 \nabla^2 f(x_t^k)(x^k - x^*) dt \right\|  \nonumber\\
   	&\leq \|x^k-x^*\| \int_0^1 \|\nabla^2 f(x_t^k)\| dt  \nonumber \\
   	&\leq M_f \| x^k-x^* \|
   	\end{align}
   and
   	\begin{align*}
   	 \|\nabla_{\Gamma_k} f({u}^k) - \nabla_{\Gamma_k} f({x}^*)\|
    &= \left\| \int_0^1 \nabla^2_{\Gamma_k, \Gamma_k} f({u}_t^k)({u}^k - {x}^*)_{\Gamma_k} dt \right\|  \nonumber\\
   	&\leq \|u^k-x^*\| \int_0^1 \|\nabla^2_{\Gamma_k, \Gamma_k} f({u}_t^k)\| dt  \nonumber \\
   	&\leq M_f \| u^k-x^* \|.
   	\end{align*}
   	This together with  \eqref{eq:Hessian_inverse_41} and \eqref{eq:dk1} gives that
   	\begin{align*}
   \|d^k\|\leq\frac{M_f}{l\iota}\|u^k-x^*\|\rightarrow0.
   \end{align*}
   	Using the second-order Taylor expansion of $f$ at $u^k$, we obtain that
   	\begin{equation}\label{eq:39:0}
   		f(v^k)=f(u^k)+\langle\nabla f(u^k),d^k\rangle +
   		\frac12(d^k)^T\nabla^2f(u^k)d^k + o(\|d^k\|^2).
   	\end{equation}
   	It follows from \eqref{hk:0} and $d^k=v^k-u^k$ that $d^k_{\overline{\Gamma}_k} = 0$. Therefore, \eqref{eq:39:0} can be rewritten as
   	\begin{equation}\label{eq:40:1}
   		\begin{aligned}
   		D_f(v^k,u^k)
   		&=f(v^k)-f(u^k)-\langle\nabla f(u^k),d^k\rangle \\
   		&=f(v^k)-f(u^k)-\langle\nabla_{\Gamma_k} f(u^k),d^k_{\Gamma_k}\rangle \\
   		&=\frac12(d^k_{\Gamma_k})^TH_kd^k_{\Gamma_k}
   		+ o(\|d^k\|^2),
   			\end{aligned}
   	\end{equation}
   	and
   	\begin{equation}
   		\label{eq:41:0}
   		\begin{aligned}
   			f(v^k)-f(u^k)
   			&=\langle \nabla_{\Gamma_k}f(u^k), d^k_{\Gamma_k} \rangle
   			+\frac12 (d^k_{\Gamma_k})^TH_kd^k_{\Gamma_k}+o(\|d^k\|^2)\\
   			&=-\frac12 (d^k_{\Gamma_k})^TH_kd^k_{\Gamma_k} + o(\|d^k\|^2),
   		\end{aligned}
   	\end{equation}
   	where the last equality follows from \eqref{hk:0}.
   	From both \eqref{eq:40:1} and \eqref{eq:41:0}, we obtain
   \begin{align*}
   	f(v^k)-f(u^k)=-D_f(v^k,u^k) + o(\|d^k\|^2).
\end{align*}
   	Since $f$ is $l$-Legendre strongly convex with respect to $h$, we have
   	\[
   	D_f(v^k,u^k) \geq lD_h(v^k,u^k),
   \]
and so
   	\begin{equation}\label{eq:43:0}
   		f(v^k)-f(u^k) \leq-lD_h(v^k,u^k) + o(\|d^k\|^2).
   	\end{equation}
   	By the $\iota$-strong convexity of $h$, we deduce that
   	$$D_h(v^k,u^k) \geq \frac{\iota}{2}\|v^k-u^k\|^2 \geq \frac{\iota}{2}\|d^k\|^2.$$
   	Due to $0<\sigma<l$, we get that for sufficiently large $k$,
   	\[
   	o(\|d^k\|^2) \leq \left(l-\sigma \right)
   	D_h(v^k,u^k).
   	\]
   	In turn, \eqref{eq:43:0} yields
   	\[
   f(v^k) \leq f(u^k) - \sigma D_h(v^k,u^k).
   \]
   	Therefore, the Newton pursuit step
   	satisfies the acceptance criterion
   	for all sufficiently large $k$.	Since it is also eventually activated, we conclude that $x^{k+1}=v^k$ for all sufficiently large $k$. Hence, (ii) is true.

   (iii)
   	We first show a local error bound
   	for the Bregman projection step. For convenience, define
   	\[
   	e_k =\|x^k-x^*\|,
   	\,\,\,
   	z_k =\|u^k-x^*\|,
   	\,\,\,
   	m_k =\|u^k-x^k\|.
   	\]
   	Taking $x=x^k$, $x_s^+=u^k$, $u=x^*$ and $\lambda=\alpha_k$ in Lemma \ref{lem:3.2}, we have
   \begin{align}\label{eq:50:1}
   	D_h(u^k,x^k)\leq D_h(x^*,x^k)+\alpha_k
   	\langle\nabla f(x^k),x^*-u^k\rangle.
   \end{align}
   From (i) and \eqref{eq:restricted_stationarity} we konw that if $\|x^*\|_0=s$, $\operatorname{supp}(u^k)=\Gamma^*$
   	and $\nabla_{\Gamma^*}f(x^*)=0$.
   	If $\|x^*\|_0<s$, $\nabla f(x^*)=0$. Therefore, for all sufficiently large $k$,
   	$$\langle\nabla f(x^*),x^*-u^k\rangle=0$$
   	and
   	\begin{align}\label{eq:gradient_bound_41}
   		\langle\nabla f(x^k),x^*-u^k\rangle
   	     &=
   		\langle
   		\nabla f(x^k)-\nabla f(x^*),
   		x^*-u^k
   		\rangle
   		\nonumber\\
   		&\leq
   		\|\nabla f(x^k)-\nabla f(x^*)\|
   		\|u^k-x^*\| \nonumber\\
   		&\leq   M_f e_k z_k,
   	\end{align}
   where the last inequality follows from \eqref{eq:so-3:00}.
   	 Then we have
   	\[
   	f(x^*)=f(x^k)+\langle \nabla f(x^k),x^*-x^k\rangle
   	+
   	\int_0^1(1-t)\,(x^*-x^k)^T\nabla^2 f(x_t)(x^*-x^k)dt,
   	\]
   	where $x_t =x^k+t(x^*-x^k), t\in[0,1]$.
   	It is equivalent to
   	\[
   	D_f(x^*,x^k)=\int_0^1(1-t)\,(x^*-x^k)^T\nabla^2 f(x_t)(x^*-x^k)\,dt.
   	\]
   Therefore, one can get that
   	\begin{align*}
   		D_f(x^*,x^k)
   		\le
   		\int_0^1(1-t)M_f\|x^*-x^k\|^2dt
   		=
   		\frac{M_f}{2}e_k^2.
   	\end{align*}
   	According to $lD_h(x^*,x^k)
   	\leq D_f(x^*,x^k)$, we have
   	\begin{equation}\label{2222}
   		D_h(x^*,x^k)\leq  \frac{1}{l}D_f(x^*,x^k)\leq \frac{M_f}{2l}e_k^2.
   	\end{equation}
   Again from the $\iota$-strong convexity of $h$, one has
   \begin{equation}\label{eq:Bregman_lower:2}
   	D_h(u^k,x^k)
   	\geq\frac{\iota}{2}m_k^2.
   \end{equation}
   Combining \eqref{eq:50:1}, \eqref{eq:gradient_bound_41}, \eqref{2222} and \eqref{eq:Bregman_lower:2}, we obtain
   \begin{align}\label{eq:projection:1}
   \frac{\iota}{2}m_k^2
    \leq
   \frac{M_f}{2l}e_k^2
   +\alpha_k M_f e_k z_k
    \leq  \frac{M_f}{2l}e_k^2
   	+\tau M_f e_k z_k,
   \end{align}
   where the last inequality results from  $\alpha_k\leq\tau$.
   	Taking into account that $z_k=\|u^k-x^*\|\leq \|u^k-x^k\|+\|x^k-x^*\|= m_k+e_k$, \eqref{eq:projection:1}
   implies
   \[
   \frac{\iota}{2}m_k^2
   \leq
   \frac{M_f}{2l}e_k^2
   +\tau M_f e_k(m_k+e_k),
   \]
  and so,
    \begin{equation*}
   m_k^2
   \leq
   \left(
   \frac{M_f}{\iota l}
   +\frac{2\tau M_f}{\iota}
   \right)e_k^2
   +
   \frac{2\tau M_f}{\iota}e_k m_k.
  \end{equation*}
   Let
   $
   A=
   \frac{M_f}{\iota l}
   +\frac{2\tau M_f}{\iota}$ and
   $B=\frac{2\tau M_f}{\iota}
   $.
   Then
   \begin{equation*}
   m_k^2-Be_km_k-Ae_k^2\leq0.
   \end{equation*}
  Moreover, we get
   \[
   m_k
   \leq
   \frac{B+\sqrt{B^2+4A}}{2}e_k.
   \]
   Since $z_k\leq e_k+m_k$,
   we conclude that
   \begin{equation}\label{eq:projection_error_bound_41}
   	\|u^k-x^*\|
   	\leq C_0\|x^k-x^*\|,
   \end{equation}
   where $	C_0=
   1+\frac{B+\sqrt{B^2+4A}}{2}>0$.

   	It follows from (i) and 	$\operatorname{supp}(u^k)\subseteq\Gamma_k$ that $	\operatorname{supp}(u^k-x^*)
   	\subseteq\Gamma_k$.
   	By (ii), we conclude that $x^{k+1}=v^k$ for all
   	sufficiently large $k$. By the Newton pursuit step of Algorithm \ref{alg:GPNP},
   	$\operatorname{supp}(v^k)\subseteq\Gamma_k$,
   	and so
   	\[
   	\|v^k-x^*\|
   	=
   	\|(v^k-x^*)_{\Gamma_k}\|.
   	\]
   	From \eqref{hk:0}, we have $v^k_{\Gamma_k}
   	=
   	u^k_{\Gamma_k}
   	-(H^k)^{-1}\nabla_{\Gamma_k}f(u^k)$. Using \eqref{eq:restricted_stationarity} yields that
   	\begin{align*}
   		(v^k-x^*)_{\Gamma_k}=
   		(u^k-x^*)_{\Gamma_k}-(H^k)^{-1}
   		\big(
   		\nabla_{\Gamma_k}f(u^k)
   		-\nabla_{\Gamma_k}f(x^*)
   		\big).
   	\end{align*}
   	Besides, we have
   	\begin{align*}
   		\nabla_{\Gamma_k}f(u^k)
   		-\nabla_{\Gamma_k}f(x^*)=
   		\int_0^1
   		\nabla^2_{\Gamma_k\Gamma_k}f(u_t^k)
   		(u^k-x^*)_{\Gamma_k}\,dt,
   	\end{align*}
   	where $u_t^k=x^*+t(u^k-x^*), t\in[0,1]$.
   	Hence,
   	\begin{align}
   		(v^k-x^*)_{\Gamma_k}
   		&=
   	(u^k-x^*)_{\Gamma_k}-(H^k)^{-1}
   		\int_0^1
   	\nabla^2_{\Gamma_k\Gamma_k}f(u_t^k)
   	(u^k-x^*)_{\Gamma_k}\,dt\nonumber\\
   			&=
   		(H^k)^{-1}
   		\int_0^1
   		\left[
   		H^k-
   		\nabla^2_{\Gamma_k\Gamma_k}f(u_t^k)
   		\right]\cdot
   		(u^k-x^*)_{\Gamma_k}\,dt.
   		\label{eq:Newton_error_identity_41}
   	\end{align}
   	By the local $s$-RHLC property,
   	for all sufficiently large $k$, we have
   	\begin{align}
   		\left\|
   		H^k-
   		\nabla^2_{\Gamma_k\Gamma_k}f(u_t^k)
   		\right\|
        &=
   		\left\|
   		\nabla^2_{\Gamma_k\Gamma_k}f(u^k)
   		-
   		\nabla^2_{\Gamma_k\Gamma_k}f(u_t^k)
   		\right\|
   		\nonumber\\
   		&\leq
   		M_s^*\|u^k-u_t^k\|
   		\nonumber\\
   		&=
   		(1-t)M_s^*\|u^k-x^*\|.
   		\label{eq:restricted_Hessian_bound_41}
   	\end{align}
   	Combining
   	\eqref{eq:Hessian_inverse_41},
   	\eqref{eq:Newton_error_identity_41}
   	and
   	\eqref{eq:restricted_Hessian_bound_41},
   	we obtain
   	\begin{align*}
   		\|v^k-x^*\|=\|(v^k-x^*)_{\Gamma_k}\|
   		&\leq
   		\frac{M_s^*}{l\iota}
   		\|u^k-x^*\|^2
   		\int_0^1(1-t)\,dt
   		\nonumber\\
   		&=
   		\frac{M_s^*}{2l\iota}
   		\|u^k-x^*\|^2.
   	\end{align*}
   Applying
   	\eqref{eq:projection_error_bound_41}
   	gives
   	\begin{align*}
   		\|x^{k+1}-x^*\|
   		=
   		\|v^k-x^*\|
   		\leq
   		\frac{M_s^*}{2l\iota}
   		\|u^k-x^*\|^2
   		\leq
   		\frac{M_s^*C_0^2}{2l\iota}
   		\|x^k-x^*\|^2,
   	\end{align*}
   		for all sufficiently large $k$.
   	Then, there exists an integer
   	$K>0$ such that
   	\[
   	\|x^{k+1}-x^*\|
   	\leq C\|x^k-x^*\|^2,
   	\,\,\, \forall k\geq K,
   	\]
   	where $C =\frac{M_s^*C_0^2}{2l\iota}>0$.
   Consequently, $\{x^k\}$ converges
   	Q-quadratically to $x^*$. The proof is completed.
   \end{proof}

   \begin{remark}
   	\begin{enumerate}
   		\item[\rm(i)]
   	Theorem \ref{th:4.1} shows that the BGPNP can safely switch from Bregman gradient projection to Newton pursuit step. It justifies the hybrid design mechanism of the BGPNP that the Bregman projection step guarantees feasibility and global convergence, while the Newton pursuit step improves the local convergence rate after support stabilization.

   \item[\rm(ii)] If $h(x)=\frac{\|x\|^2}{2}$, the conditions of Theorem \ref{th:4.1} reduce to that of  \cite[Theorem 3.2]{zhou2022gradient}. So, \cite[Theorem 3.2]{zhou2022gradient} is a special case of Theorem \ref{th:4.1}. In other words, Theorem \ref{th:4.1} establishes the local quadratic convergence of the proposed projection gradient-Newton-type algorithm under flexible conditions. Besides, the proof of Theorem \ref{th:4.1} (ii) is different from that of \cite[Theorem 3.2(ii)]{zhou2022gradient} even if $h(x)=\frac{\|x\|^2}{2}$ because the proof of Theorem \ref{th:4.1} (ii) is obtained by the definition of the stationary point.
   	\end{enumerate}
   	\end{remark}

\section{Numerical Experiments}
\label{sec:numerical_experiments}

In this section, we conduct numerical experiments to demonstrate the
performance of BGPNP for sparsity-constrained compressed sensing problems.
The initial point is set to $x^{0}=0$. Unless otherwise specified, the main
parameters of BGPNP are chosen as follows: $\tau=5, \gamma=0.5, \sigma=10^{-4}, \kappa=10^{-2}$ and $\epsilon= 10^{-5}$.
All algorithms are implemented in MATLAB $2021a$ on a laptop equipped with an Intel(R) Core(TM) i5-8250U processor running at 1.60 GHz.

We compare BGPNP with GPNP\cite{zhou2022gradient} and several representative greedy
methods, including
GP\cite{BlumensathDavies2008GP}, NIHT\cite{blumensath2010normalized}, OMP\cite{tropp2007omp},
CoSaMP\cite{needell2009cosamp}, HTP\cite{foucart2011htp}, and SP\cite{dai2009subspace}.
Numerical experiments validate the competitive performance of BGPNP
compared with existing leading solvers.

\subsection{Solving  Compressed Sensing Problems}
\label{subsec:solving_cs}

The compressed sensing (CS) problem can be formulated as problem \eqref{qu:1}
with the objective function
\begin{align*}
	f(x)=f_{\rm cs}(x)=
	\frac{1}{2}\|Ax-b\|^{2},
	\label{eq:cs_objective}
\end{align*}
where $A\in\mathbb{R}^{m\times n}$ is the sensing matrix and
$b\in\mathbb{R}^{m}$ is the observation vector.

We apply  BGPNP to solve the CS problem directly rather than its approximate problem such as $\ell_{1}$-norm regularization problem.
For BGPNP, we employ the following Bregman kernel function
\begin{equation}
	h_{\beta}(x)
	=
	\frac{1}{2}\|x\|^{2}
	+
	\frac{2\beta}{3}\|x\|_{3/2}^{3/2},
	\,\, \beta\geq 0.
	\label{eq:bregman_kernel}
\end{equation}
Then
 $	\nabla h_{\beta}(x)= x+
	\beta\,\operatorname{sign}(x)
	\odot\sqrt{|x|}$,
where $\odot$ denotes componentwise multiplication. Moreover, the inverse
mapping of $\nabla h_{\beta}$ admits the explicit form
\begin{equation*}
	[\nabla h_{\beta}]^{-1}(y)
	=
	\operatorname{sign}(y)\odot
	\left(
	\frac{\sqrt{\beta^{2}+4|y|}-\beta}{2}
	\right)^{2}.
	\label{eq:bregman_inverse}
\end{equation*}
Therefore, the corresponding Bregman gradient step can be computed explicitly
without solving an additional nonlinear subproblem. In particular, if
$\beta=0$,  the kernel function \eqref{eq:bregman_kernel}  reduces to the function $h_{0}(x)=\frac{1}{2}\|x\|^{2}$.
Then, BGPNP recovers to GPNP\cite{zhou2022gradient}. The Bregman parameter $\beta$ is specified separately in each experiment.   The numerical influence of
$\beta$ will be investigated in the next subsection.

We consider two types of CS problems: randomly generated Gaussian data and
two-dimensional image data. For the former, we consider exact recovery
\[
b=Ax^{*},
\]
while for the image data, we consider noisy recovery
\[
b=Ax^{*}+\chi,
\]
where $\chi$ denotes measurement noise.
\vskip2mm

\subsubsection{Testing Examples}
\label{subsubsec:testing_examples}

\begin{example}{\rm (Gaussian matrix)}{\rm \cite{zhou2022gradient}}\label{exam:4.1}
Let $A\in\mathbb{R}^{m\times n}$ be a Gaussian sensing matrix whose entries
are independently and identically distributed according to
$\mathcal{N}(0,1)$. Each column of $A$ is then normalized to have unit
Euclidean norm. The support of the ground-truth signal $x^{*}$ is selected
randomly, and its $s$ nonzero entries are independently generated from
$\mathcal{N}(0,1)$. Finally, the observation vector is given by
$b=Ax^{*}$.
\end{example}

\begin{example}{\rm (Image data)}{\rm \cite{zhou2022gradient}}\label{exam:4.2}
Images are generally not sparse in the spatial domain but can possess sparse
representations under wavelet transforms. We use a wavelet transform
$W(\cdot)$ and denote the corresponding sparse coefficients by $x^{*}$.
The sensing operator takes the form
   $A=FW^{-1}$,
where $F$ denotes a partial Fourier sampling operator and $W^{-1}$ denotes
the inverse wavelet transform. The noisy observation is
\[
b=Ax^{*}+\zeta,
\]
where the entries of $\zeta$ are generated according to $\zeta_i\sim nf\,\mathcal{N}(0,1)$. We consider two noise factors, $nf=0.05$ and $nf=0.1$. A $256\times256$ gray-scale image and a $256\times256$ colour image are considered in the experiments.
\end{example}

For Example \ref{exam:4.1}, the recovery accuracy is measured by the relative error
\begin{equation*}\label{eq:relative_error}
	{\rm ReEr}
	=
	\frac{\|x-x^{*}\|}
	{\|x^{*}\|}.
\end{equation*}
A recovery is regarded as successful iff
	${\rm ReEr}<10^{-2} $.

For Example \ref{exam:4.2}, the reconstruction quality is evaluated by the peak
signal-to-noise ratio
\begin{equation*}
	{\rm PSNR}
	=
	10\log_{10}
	\left(
	\frac{N}
	{\|I-\widehat{I}\|_{F}^{2}}
	\right),
	\label{eq:PSNR}
\end{equation*}
where $I$ and $\widehat{I}$ denote the original and reconstructed images,
respectively, and $N$ denotes the total number of pixels.

\subsubsection{Effect of Parameter $\beta$}
\label{subsubsec:beta_effect}

We now consider the effect of the Bregman parameter $\beta$ on the
performance of  Algorithm \ref{alg:GPNP} by
running $100$ independent trials with fixed $n=256, m=64$ and
recording the corresponding success rates
(which is defined as the percentage of the number of successful recoveries over all trials) at sparsity levels
$s$ from $5$ to $35$.
  The following values $\beta\in\{0.005,\ 0.01,\ 0.05,\ 0.1,\ 0.2,\ 0.5,\ 1\}$ are tested.

\begin{figure}[htbp]
	\centering
	\includegraphics[width=0.8\textwidth]{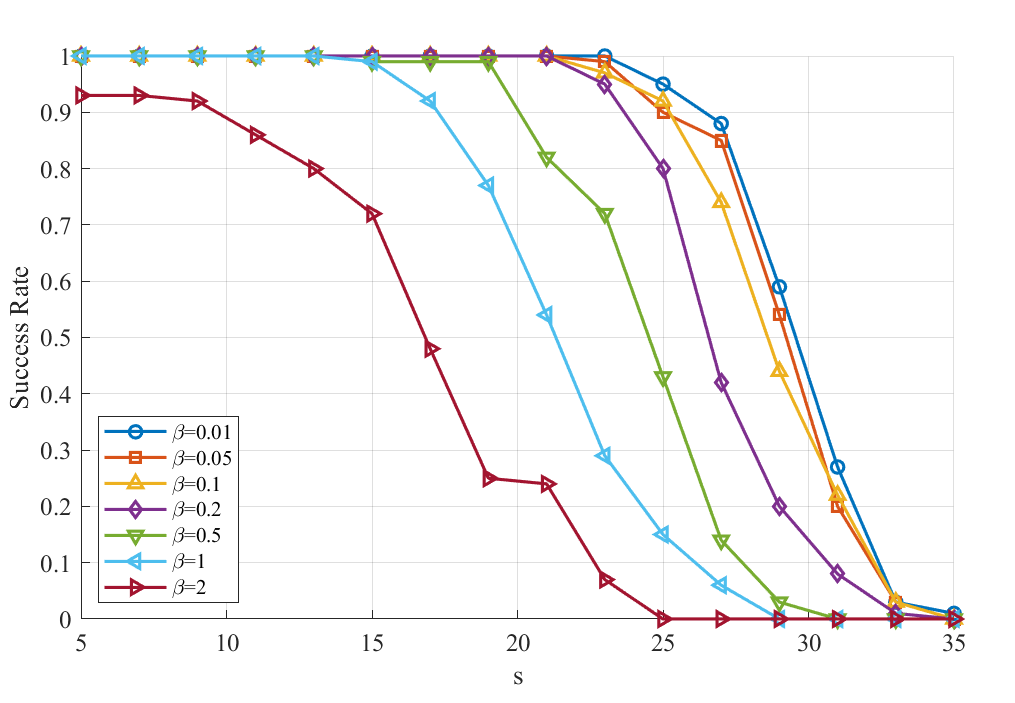}
	\caption{Recovery success rate of BGPNP with different parameter
		$\beta$ versus the sparsity level $s$.}
	\label{fig:beta_success}
\end{figure}

Figure~\ref{fig:beta_success} reports the recovery success rates. For small
values of $s$, all tested values of $\beta$ achieve nearly $100\%$ successful
recovery. As $s$ increases, the recovery becomes more difficult and the
difference among different values of $\beta$ becomes evident. Relatively
small values of $\beta$ generally preserve better recovery performance,
whereas large values such as $\beta=1$ and $\beta=2$ lead to an earlier
decrease in the success rate.

Table~\ref{tab:beta_effect} presents the relative errors and
CPU times for different parameter $\beta$. For relatively small
values of $\beta$, BGPNP achieves highly accurate recovery for a wide range
of sparsity levels. In contrast, excessively large values of $\beta$ may
deteriorate the recovery accuracy when the problem becomes more difficult.
The results suggest that small positive values such as
$\beta=0.005, 0.01, 0.05$, and $0.1$ provide suitable choices for the subsequent
comparisons.

\begin{table}[htbp]
	\centering
	\caption{Effect of the Bregman parameters $\beta$ on the relative error and
		CPU time of BGPNP.}
	\label{tab:beta_effect}
	
\resizebox{\textwidth}{!}{
\begin{tabular}{c ccccccc}
	\toprule
	& \multicolumn{7}{c}{ReEr} \\
	\cmidrule(lr){2-8}
	$\beta\backslash s$
	& 5 & 10 & 15 & 20 & 25 & 30 & 35 \\
	\midrule
	
0.005
& 6.03e-17 & 1.29e-16 & 2.34e-16 & 3.14e-16
& 3.55e-16 & 4.78e-16 & 1.17e-15 \\

0.01
& 6.03e-17 & 1.29e-16 & 2.15e-16 & 3.11e-16
& 3.52e-16 & 5.41e-16 & 0.1339 \\

0.05
& 6.03e-17 & 1.29e-16 & 2.34e-16 & 2.94e-16
& 3.52e-16 & 5.85e-16 & 1.17e-15 \\

0.1
& 6.03e-17 & 1.29e-16 & 2.46e-16 & 3.14e-16
& 3.55e-16 & 5.63e-16 & 0.3045 \\

0.2
& 6.03e-17 & 1.29e-16 & 2.46e-16 & 2.97e-16
& 3.55e-16 & 6.01e-16 & 0.0004 \\

0.5
& 6.03e-17 & 1.29e-16 & 2.46e-16 & 2.94e-16
& 3.55e-16 & 7.34e-16 & 0.3517 \\

1
& 6.03e-17 & 1.76e-16 & 2.15e-16 & 3.48e-16
& 4.59e-16 & 0.3078 & 0.2862 \\

\midrule
& \multicolumn{7}{c}{Time (seconds)} \\
\cmidrule(lr){2-8}
$\beta\backslash s$
	& 5 & 10 & 15 & 20 & 25 & 30 & 35 \\
	\midrule
	
	0.005
	& 0.0048 & 0.0101 & 0.0129 & 0.0265
	& 0.1009 & 0.2838 & 0.1671 \\
	
	0.01
	& 0.0088 & 0.0072 & 0.0141 & 0.0314
	& 0.2031 & 0.0689 & 3.4546 \\
	
	0.05
	& 0.0033 & 0.0046 & 0.0086 & 0.0924
	& 0.8522 & 0.1069 & 0.1387 \\
	
	0.1
	& 0.0028 & 0.0063 & 0.0077 & 0.0081 & 0.0619 & 0.0967 & 2.1964 \\
	
	0.2
	& 0.0042 & 0.0056 & 0.0092 & 0.0092 & 0.0265 & 0.0306 & 0.0420 \\
	
	0.5
	& 0.0043 & 0.0038 & 0.0053 & 0.0242
	& 0.0177 & 0.0177 & 0.0171 \\
	
	1
	& 0.0019 & 0.0038 & 0.0067 & 0.0131
	& 0.0174 & 1.5273 & 0.0069 \\
	\bottomrule
\end{tabular} }

\end{table}

\subsubsection{Numerical Comparisons with Benchmark Methods}
\label{subsubsec:benchmark}

In this subsection,
we compare BGPNP with seven representative sparse recovery methods:
GP\cite{BlumensathDavies2008GP}, NIHT\cite{blumensath2010normalized}, OMP\cite{tropp2007omp}, CoSaMP\cite{needell2009cosamp}, HTP\cite{foucart2011htp}, SP\cite{dai2009subspace},
 and GPNP\cite{zhou2022gradient}. For BGPNP, the parameter $\beta=0.005, 0.01, 0.05, 0.1, 0.5$
are considered in the success-rate experiments. For each trial, all methods
are applied to the same sensing matrix, observation vector, and ground-truth
signal to ensure a fair comparison.

%
\noindent
\textbf{a) Recovery Success Rate Versus Sparsity Level}

\begin{figure}[htbp]
	\centering
	\includegraphics[width=0.8\textwidth]{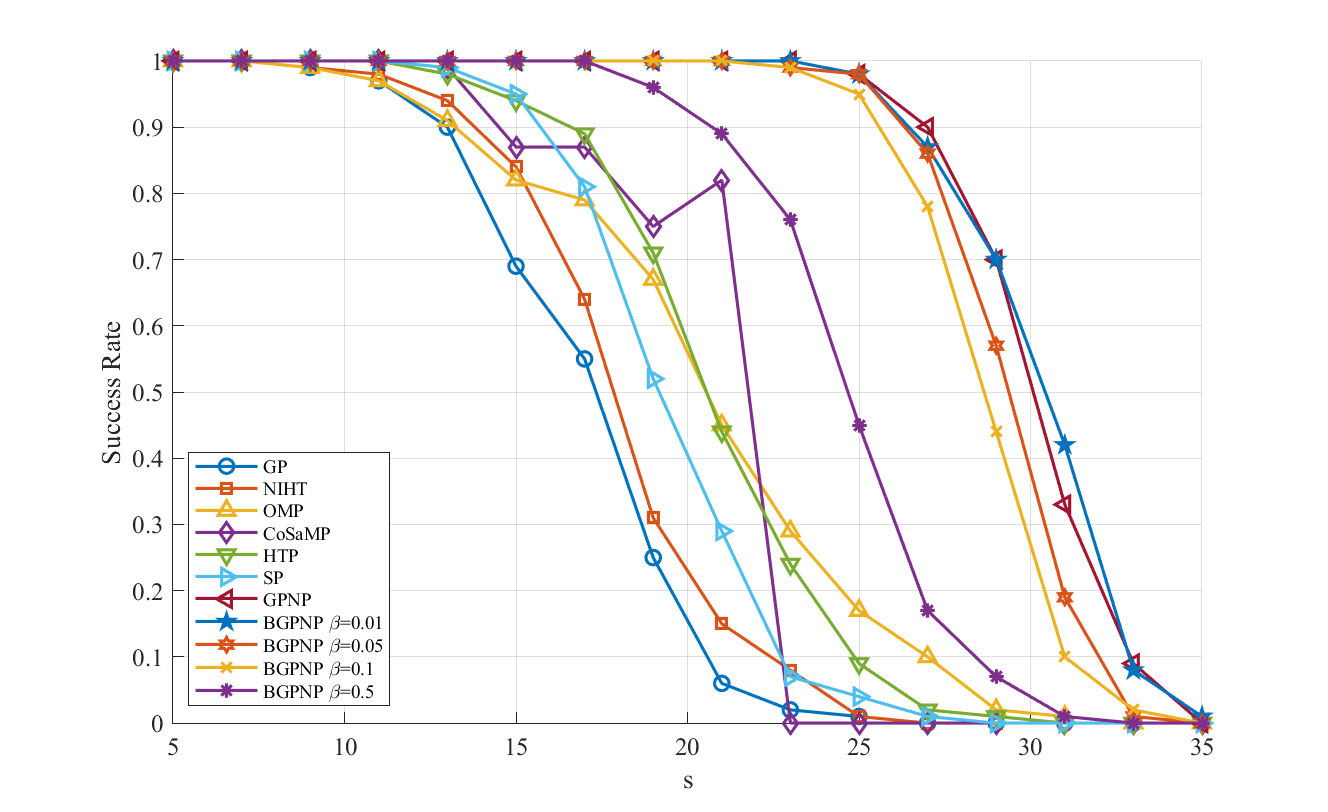}
	\caption{Recovery success rate versus the sparsity level $s$ for
		Example~\ref{exam:4.1}}
	\label{fig:success_s}
\end{figure}

We first fix $n=256, m=64$, and vary the sparsity level from $s=5$ to $s=35$. The corresponding recovery success rates are reported in Figure~\ref{fig:success_s}.

As expected, increasing $s$ makes the recovery problem more difficult.
Several classical greedy methods exhibit an earlier decrease in recovery
success rate. In comparison, GPNP and BGPNP with suitable values of $\beta$
maintain relatively high success rates over a wider range of sparsity levels.
In particular, the BGPNP variants with small positive $\beta$ show
competitive recovery performance in the transition region.

\medskip
\noindent
\textbf{b) Recovery Success Rate Versus Sampling Ratio.}

\begin{figure}[htbp]
	\centering
	\includegraphics[width=0.8\textwidth]{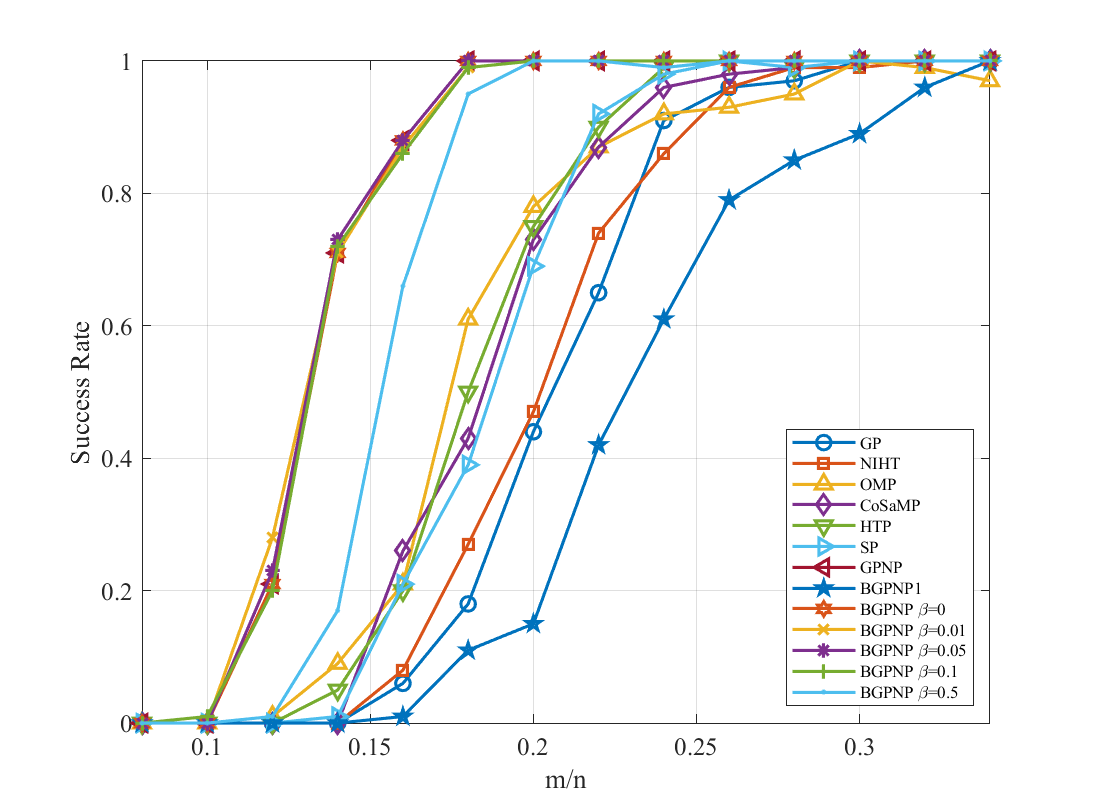}
	\caption{Recovery success rate versus the sampling ratio $m/n$ for
		Example~\ref{exam:4.1}}
	\label{fig:success_ratio}
\end{figure}

We next fix $n=256, s=13$, and vary the sampling ratio $m/n$ from $0.08$ to $0.34$.
Figure~\ref{fig:success_ratio} presents the corresponding recovery success
rates.

All methods show increasing success rates as the number of measurements
increases. The differences are mainly observed in the transition region,
where GPNP and several BGPNP variants achieve high recovery probabilities
with relatively fewer measurements. When the sampling ratio becomes
sufficiently large, almost all tested algorithms achieve successful
recovery.

\medskip
\noindent
\textbf{c) Higher-Dimensional Problems.}

We further validate the performance of the above algorithms for
higher-dimensional CS problems. The signal dimension is chosen from $n\in\{5000, 10000$, $15000\}$, and two sparsity levels $s=0.25\%n$
and $s=5\%n$ are considered.

Table~\ref{tab:large025} reports the numerical results for the sparsity level $s=0.25\%n$.
OMP, CoSaMP, HTP, SP, GPNP, and BGPNP achieve highly accurate
recovery. In terms of CPU time, BGPNP performs better than GP, NIHT, OMP, CoSaMP, HTP, SP and GPNP.
In particular, BGPNP with $\beta=0.05$ performs the best among these algorithms from the CPU time. Concretely, BGPNP with $\beta=0.05$ requires
$0.0620$, $0.2303$, and $0.4886$ seconds for
$n=5000$, $10000$, and $15000$, respectively, compared with
$0.0830$, $0.3626$, and $0.7004$ seconds required by GPNP.

%
%
%
%
%
%
%
%
%
%
%
%
%

\begin{table}[htbp]
	\centering
	\caption{Relative errors and CPU time for Example~\ref{exam:4.1} with
		$s=0.25\%n$.}
	\label{tab:large025}
	\resizebox{\textwidth}{!}{
	\begin{tabular}{l ccc ccc}
		\toprule
		& \multicolumn{3}{c}{ReEr}
		& \multicolumn{3}{c}{Time (seconds)} \\
		\cmidrule(lr){2-4}
		\cmidrule(lr){5-7}
		
		Algs.$\backslash n$
		& 5000 & 10000 & 15000
		& 5000 & 10000 & 15000 \\
		\midrule
		
		GP
		& 5.85e-12 & 3.80e-12 & 3.07e-12
		& 1.1943 & 5.1674 & 11.992 \\
		
		NIHT
		& 2.16e-09 & 1.42e-09 & 1.02e-09
		& 0.2730 & 1.1566 & 2.6669 \\
		
		OMP
		& 3.79e-16 & 6.14e-16 & 7.63e-16
		& 0.2768 & 2.3167 & 7.3489 \\
		
		CoSaMP
		& 3.78e-16 & 6.89e-16 & 8.47e-16
		& 0.0918 & 0.3703 & 0.8696 \\
		
		HTP
		& 3.90e-16 & 6.36e-16 & 8.16e-16
		& 0.0758 & 0.3232 & 0.7519 \\
		
		SP
		& 4.31e-16 & 7.04e-16 & 7.88e-16
		& 0.0765 & 0.3200 & 0.8022 \\
		
	    AIHT
	    & 9.59e-11 & 4.98e-12 & 8.78e-12
	 &0.4573 & 1.4796 & 7.2874 \\
		
		ALPS
		& 1.51e-08  & 9.50e-09 & 2.64e-08 & 0.5982 & 2.3761 & 12.989\\
		
		FHTP
		& 2.31e-11  & 6.95e-11  & 6.75e-11 & 0.4003 & 1.6514 & 5.6514\\
		
		GPNP
		& 3.44e-16 & 4.99e-16 & 3.55e-16
		& 0.0830 & 0.3626 & 0.7004 \\
		
		YALL1    & 2.63e-07 & 2.69e-07 & 2.69e-07 & 4.9269 & 25.099 & 89.809 \\
		L1Magic  & 3.74e-06 & 1.75e-06 & 5.06e-06 & 24.822 & 218.20 & 1019.5 \\
		SPGL1    & 1.99e-07 & 5.31e-08 & 1.10e-07 & 0.6589 & 2.0468 & 6.5718 \\
		GPSR     & 7.91e-03 & 9.84e-03 & 1.07e-02 & 1.0543 & 4.9150 & 13.641 \\
		IRL1     & 8.74e-06 & 7.02e-06 & 1.07e-06 & 4.4174 & 17.310 & 74.263 \\
		DCL1L2   & 1.52e-06 & 7.64e-07 & 6.91e-07 & 82.667 & 332.17 & 805.86 \\
		IRLSLq   & 1.39e-04 & 2.50e-03 & 2.38e-02 & 6.9445 & 38.617 & 907.65 \\
		AMP      & 7.80e-09 & 5.70e-09 & 5.92e-09 & 0.4726 & 1.3913 & 3.6205 \\
		SAMP     & 3.46e-16 & 5.79e-16 & 9.04e-16 & 0.3805 & 2.1930 & 8.1579 \\
		
			BGPNP($\beta=0.005$)
		& 3.04e-16 & 4.53e-16 & 2.67e-16
		& 0.0381 & 0.2327 & 1.4434 \\
		
			BGPNP($\beta=0.01$)
		& 3.44e-16 & 4.99e-16 & 3.55e-16
		& 0.0747 & 0.2619 & 0.5681 \\
		
		BGPNP($\beta=0.05$)
		& 3.44e-16 & 4.98e-16 & 3.55e-16
		& 0.0620 & 0.2303 & 0.4886 \\
		
		BGPNP($\beta=0.1$)
		& 3.44e-16 & 4.98e-16 & 3.56e-16
		& 0.0671 & 0.2439 & 0.5128 \\

		\bottomrule
	\end{tabular} }
	
\end{table}

Table~\ref{tab:large5} reports the numerical results for the sparsity level
$s=5\%n$. The computational costs of several classical methods increase
rapidly with the dimension. For $n=15000$, GP and OMP require
$184.16$ and $693.77$ seconds, respectively, while GPNP requires
$6.7427$ seconds with a relative error
$1.41\times10^{-15}$. In contrast, the CPU time of BGPNP with $\beta=0.01, 0.05$ and $0.1$ are fewer than that of GP, NIHT, OMP, CoSaMP, HTP, SP and GPNP for $n=10000$ and $15000$. In particular,  for $n=10000$,
 BGPNP with $\beta=0.05$ requires only
$0.9718$ seconds while maintaining a relative error
$1.18\times10^{-15}$; For $n=15000$,
 BGPNP with $\beta=0.01$ requires only
$3.0142$ seconds while maintaining a relative error
$1.40\times10^{-15}$. These demonstrate the computational advantage of
BGPNP for large-scale sparse recovery.

It is also observed that larger value  $\beta$ may affect the final
accuracy in more difficult problems. For $n=15000$, the relative errors
corresponding to $\beta=0.05$ and $\beta=0.1$ increase to
$7.23\times10^{-8}$ and $1.66\times10^{-7}$, respectively. Therefore,
a relatively small value $\beta$ provides a better balance between
accuracy and computational efficiency in this setting.

%
%
%
%
%
%
%
%
%
%
%
%

\begin{table}[htbp]
	\centering
	\caption{Relative errors and CPU time for Example~\ref{exam:4.1} with
		$s=5\%n$.}
	\label{tab:large5}
	\resizebox{\textwidth}{!}{
	
	\begin{tabular}{l ccc ccc}
		\toprule
		& \multicolumn{3}{c}{ReEr}
		& \multicolumn{3}{c}{Time (seconds)} \\
		\cmidrule(lr){2-4}
		\cmidrule(lr){5-7}
		
		Algs.$\backslash n$
		& 5000 & 10000 & 15000
		& 5000 & 10000 & 15000 \\
		\midrule
		
		GP
		& 2.97e-12 & 2.12e-12 & 1.73e-12
		& 24.221 & 82.666 & 184.16 \\
		
		NIHT
		& 9.25e-10 & 6.52e-10 & 5.33e-10
		& 1.1863 & 3.5615 & 7.6013 \\
		
		OMP
		& 7.53e-16 & 1.41e-15 & 1.67e-15
		& 28.134 & 105.84 & 693.77 \\
		
		CoSaMP
		& 7.64e-16 & 1.71e-15 & 1.99e-15
		& 5.5950 & 9.5478 & 53.338 \\
		
		HTP
		& 7.72e-16 & 1.44e-15 & 1.73e-15
		& 1.1622 & 2.4588 & 34.163 \\
		
		SP
		& 7.39e-16 & 1.51e-15 & 1.83e-15
		& 2.4775 & 4.7629 & 36.849 \\
		
		AIHT
		& 1.06e-10 & 5.04e-11  & 6.29e-11
		& 1.0185 & 3.2047 & 7.3868 \\
		
		ALPS
		& 2.53e-08 & 2.82e-08 & 2.30e-08 & 2.2384 & 8.0051 & 14.725\\
		
		FHTP
		& 8.43e-11 & 4.95e-11 & 8.22e-11 & 1.3602 & 4.8138 & 10.423\\
		
		GPNP
		& 8.55e-16 & 1.18e-15 & 1.41e-15
		& 0.8720 & 1.4734 & 6.7427 \\

		YALL1    & 3.72e-07 & 1.32e-04 &5.61e-05 & 5.6531 & 27.839 & 106.54 \\
		L1Magic  & 8.23e-06 & 4.82e-06 &1.60e-05 & 28.952 & 277.90 & 1504.0 \\
		SPGL1    & 7.61e-08 & 4.23e-08 &5.49e-08 & 1.0702 & 4.6596 & 15.562 \\
		GPSR     & 1.68e-02 & 2.01e-02 &1.94e-02 & 2.0780 & 8.5940 & 25.102 \\
		IRL1     & 6.54e-06 & 8.23e-06 &9.20e-06 & 8.0577 & 115.01 & 764.77 \\
		DCL1L2   & 5.97e-07 & 4.75e-06 &5.26e-06 & 105.94 & 611.87 & 25967  \\
		IRLSLq   & 5.12e-02 & 6.62e-02 &8.78e-02 & 98.199 & 462.40 & 1024.0 \\
		AMP      & 1.06e-08 & 1.22e-08 &1.46e-08 & 0.6438 & 3.1068 & 6.3283 \\
		SAMP     & 7.08e-16 & 1.50e-02 &1.53e-01 & 72.634 & 771.37 & 1118.7 \\

			BGPNP($\beta=0.005$)
		& 7.67e-16 & 1.10e-15 & 1.29e-15
		& 0.4742 & 1.1931 & 3.3674 \\
		
		BGPNP($\beta=0.01$)
		& 8.49e-16 & 1.18e-15 & 1.40e-15
		& 1.3145 & 1.0804 & 3.0142 \\
		
		BGPNP($\beta=0.05$)
		& 8.55e-16 & 1.18e-15 & 7.23e-08
		& 1.3578 & 0.9718 & 3.0969 \\
		
		BGPNP($\beta=0.1$)
		& 8.60e-16 & 2.46e-07 & 1.66e-07
		& 1.0780 & 1.0820 & 3.6899 \\
		
		\bottomrule
	\end{tabular}  }
	
\end{table}

\medskip
\noindent
\textbf{d) Comparisons for Solving Example~\ref{exam:4.2}.}

We now compare the above mentioned algorithms for noisy image recovery.
The sparsity levels are chosen as $s=1000, 1500, 2000$, and the noise factors $nf=0.05$ and $nf=0.1$.
Besides,
OMP and HTP are not suitable to solve Example \ref{exam:4.2} because measurement $A$ is an operator instead of a matrix; see \cite{zhou2022gradient}. So, we
exclude OMP and HTP in the  experiments of this subsection.
 The PSNR values and CPU time are reported in
Table~\ref{tab:image_results}.

For the gray-scale image with $nf=0.05$, BGPNP and GPNP generally achieve
higher reconstruction quality than SP and CoSaMP while requiring
substantially less CPU time than  that of GP and NIHT. For
$s=1500, 2000$, the values of PSNR for GPNP are $21.88$ dB and $21.66$ dB, while the corresponding values of PSNR for BGPNP with
$\beta=0.1$  are $21.90$ dB and $21.69$ dB, respectively.

For the colour image, BGPNP also shows competitive reconstruction
performance. In particular, when $nf=0.1$, BGPNP with $\beta=0.05$
achieves PSNR values of $27.41$ dB, $25.80$ dB, and $23.89$ dB for
$s=1000$, $1500$, and $2000$, respectively. These values are higher than
the corresponding PSNR values  $21.42$ dB, $20.00$ dB, and $19.11$ dB of GPNP.
Meanwhile, the CPU time of BGPNP with $\beta=0.05$ is also fewer than that of GPNP for $s=1000$, $1500$, and $2000$.
The numerical results indicate that an appropriate Bregman parameter can improve
the robustness of BGPNP for certain noisy image recovery problems.

\begin{table}[htbp]
	\centering
	\caption{PSNR values and CPU time for gray-scale and colour
		image recovery under different noise levels.}
	\label{tab:image_results}
	
	\resizebox{\textwidth}{!}{
		\begin{tabular}{l ccc ccc ccc ccc}
			\toprule
			
			&
			\multicolumn{6}{c}{Gray image}
			&
			\multicolumn{6}{c}{Colour image}
			\\
			
			\cmidrule(lr){2-7}
			\cmidrule(lr){8-13}
			
			&
			\multicolumn{3}{c}{PSNR}
			&
			\multicolumn{3}{c}{Time (seconds)}
			&
			\multicolumn{3}{c}{PSNR}
			&
			\multicolumn{3}{c}{Time (seconds)}
			\\
			
			\cmidrule(lr){2-4}
			\cmidrule(lr){5-7}
			\cmidrule(lr){8-10}
			\cmidrule(lr){11-13}
			
			$s$
			& 1000 & 1500 & 2000
			& 1000 & 1500 & 2000
			& 1000 & 1500 & 2000
			& 1000 & 1500 & 2000
			\\
			
			\midrule
			\multicolumn{13}{l}{$nf=0.05$}\\
			\midrule
			
			GP
			& 21.64 & 21.73 & 21.37
			& 31.84 & 32.43 & 33.09
			& 29.25 & 29.17 & 28.00
			& 94.68 & 104.50 & 99.47 \\
			
			NIHT
			& 21.60 & 21.61 & 20.92
			& 13.98 & 14.28 & 17.02
			& 29.00 & 28.55 & 27.43
			& 24.60 & 48.04  & 64.01  \\
			
			SP
			& 21.28 & 20.42 & 19.75
			& 2.54 & 6.84 & 14.12
			& 27.00 & 25.19 & 24.58
			& 19.80 & 53.82  & 48.79 \\
			
			CoSaMP
			& 21.02 & 20.15 & 19.84
			& 4.26 & 16.76 & 18.76
			& 26.61 & 25.52 & 24.29
			& 28.92 & 39.73  & 83.03 \\
			
		    AIHT
			& 21.67 & 21.70 & 21.34
			& 5.53 & 5.19 & 7.46
			& 29.21 & 29.11 & 27.87
			& 13.62 & 16.24  & 17.46 \\
			
			ALPS
			& 21.71 & 21.66 & 21.33
			& 3.09 & 5.83 & 5.68
			& 29.30 & 29.11 & 27.87
			& 10.58 & 12.91  & 16.73\\
			
			FHTP
			& 21.67 & 21.73 & 21.45
			& 4.44 & 8.88 & 9.41
			& 16.32 & 29.14  & 33.81
			& 16.32 & 24.16  & 33.81 \\
			
			GPNP
			& 21.65 & 21.88 & 21.66
			& 2.58 & 2.56 & 4.82
			& 29.40 & 29.63 & 28.39
			& 5.97 & 9.04 & 11.07 \\
			
			BGPNP($\beta=0.005$)
			& 21.65 & 21.87 & 21.64
			& 1.64 & 2.41 & 5.40
			& 29.40 & 22.02 & 28.42
			& 5.55 & 11.08 & 9.65 \\
			
			BGPNP($\beta=0.01$)
			& 21.65 & 21.86 & 21.67
			& 1.58 & 2.17 & 3.36
			& 29.40 & 29.71 & 21.83
			& 5.93 & 8.32  & 11.15 \\
			
			BGPNP($\beta=0.05$)
			& 21.64 & 21.86 & 21.68
			& 2.01 & 2.26 & 3.93
			& 29.40 & 29.13 & 28.23
			& 6.27 & 8.78  & 9.55 \\
			
			BGPNP($\beta=0.1$)
			& 21.64 & 21.90 & 21.69
			& 0.89 & 1.66 & 1.22
			& 29.36 & 28.72 & 27.84
			& 6.33 & 6.33  & 7.30 \\
			
			BGPNP($\beta=0.5$)
			& 21.61 & 21.73 & 21.42
			& 2.43 & 2.85 & 4.38
			& 28.76 & 28.52 & 27.57
			& 4.26 & 5.19 & 7.24 \\

			\midrule
			\multicolumn{13}{l}{$nf=0.10$}\\
			\midrule
			
			GP
			& 21.28 & 20.73 & 20.08
			& 32.75 & 31.85 & 32.73
			& 26.90 & 24.81 & 22.96
			& 102.86 & 105.61 & 106.23 \\
			
			NIHT
			& 21.13 & 20.29 & 19.57
			& 12.51 & 29.78 & 24.40
			& 25.95 & 23.85 & 22.41
			& 34.14  & 42.88  & 54.58 \\
			
			SP
			& 19.22 & 17.55 & 17.10
			& 4.57  & 26.59 & 12.68
			& 22.60 & 20.32 & 19.62
			& 17.25  & 38.64  & 50.91 \\
			
			CoSaMP
			& 19.42 & 18.13 & 17.29
			& 4.21  & 7.98  & 13.85
			& 22.42 & 20.32 & 19.23
			& 15.20  & 37.71  & 98.71\\
			
			AIHT
			& 21.28 & 20.73 & 19.86
			& 3.50  & 4.86  & 6.64
			& 26.78 & 24.61 & 22.94
			& 14.97  & 19.79  & 21.60 \\
			
			ALPS
			& 21.27 & 20.74 & 19.74
			& 2.74  & 4.00  & 5.69
			& 26.89 & 24.67 & 22.75
			& 10.48  & 14.96  & 21.08 \\
			
			FHTP
			& 21.30 & 20.73 & 20.08
			& 4.15  & 7.00  & 9.54
			& 26.89 & 24.93 & 23.10
			& 16.37  & 23.76  & 28.76 \\
			
			GPNP
			& 21.34 & 20.94 & 20.44
			& 1.75 & 2.02 & 4.72
			& 27.42 & 25.81 & 23.95
			& 6.03 & 8.80 & 14.05 \\
			
			BGPNP($\beta=0.005$)
			& 21.34 & 20.95 & 20.54
			& 1.27 & 2.64 & 2.78
			& 27.45 & 25.81 & 24.03
			& 5.69 & 8.37 & 13.23 \\
			
			BGPNP($\beta=0.01$)
			& 21.34 & 21.00 & 20.54
			& 1.30 & 1.77 & 2.80
			& 27.44 & 25.83 & 23.44
			& 5.44 & 9.57 & 15.49 \\
			
			BGPNP($\beta=0.05$)
			& 21.34 & 20.99 & 20.50
			& 2.34 & 2.76 & 3.43
			& 27.49 & 24.80 & 22.23
			& 6.16 & 8.78 & 9.26 \\
			
			BGPNP($\beta=0.1$)
			& 21.33 & 21.11 & 19.23
			& 1.92 & 2.10& 3.48
			& 27.54 & 24.12 & 22.61
			& 6.17 & 8.12 & 7.31 \\
			
			BGPNP($\beta=0.5$)
			& 21.31 & 20.62 & 20.09
			& 0.87 & 1.39 & 1.26
			& 25.72 & 24.04 & 22.50
			& 3.36 & 3.92 & 4.87 \\
			
			\bottomrule
		\end{tabular}
	}
	
\end{table}



	

\section{Conclusions}
  The BGPNP is proposed to solve sparse constrained optimization problem via the Bregman distance. By combining Bregman gradient projection with Newton pursuit, the proposed method inherits the global convergence property of first-order methods. The locally linear convergence of BGPNP without the Newton pursuit step is established under suitable conditions. Further,  when the support set becomes stable or the gradient norm is sufficiently small, the convergence rate of BGPNP can be enhanced to the locally
quadratical convergence by adding a Newton pursuit step.
  This method provides an effective algorithmic framework for nonconvex sparse optimization in non-Euclidean settings. Numerical experiments validate the competitive performance of BGPNP compared with existing leading solvers.  Further research, it is also interesting to consider the globally linear convergence and quadratic convergence of BGPNP. It also remains an open question whether linear convergence of BGPNP  can be established beyond  the convexity assumption and $\|x^{*}\|_{0}=s$.
\vskip3mm

\noindent{\bf Acknowledgements}.
The initial manuscript was finalized in December 2025.
This paper was partially supported by
the National Natural Science Foundation of China (Nos. 12471299, 12071379), the Fundamental Research Funds for the Central Universities (Nos. SWU-KF26001, SWU-KT25018),
 the Natural Science Foundation of Chongqing (No. CSTB2025NSCQ-GPX0529) and Research Project on Education and Teaching Reform of Southwest University  (2025JY004).
 This work was also partially supported by National Natural Science Foundation of China (No. 12271061).


\end{document}